# Mordukhovich derivatives of the metric projection operator in uniformly convex and uniformly smooth Banach spaces

Jinlu Li

Department of Mathematics
Shawnee State University
Portsmouth, Ohio 45662
USA

**Abstract**

In this paper, we investigate the properties of the Mordukhovich derivatives of the metric projection operator onto closed balls, closed and convex cylinders and positive cones in uniformly convex and uniformly smooth Banach spaces. We find the exact expressions for Mordukhovich derivatives of the metric projection operator.

1. Introduction

In the operator theory, it is clear that the metric projection operator is one of the most important operators, which has been studied by many authors with a long history (see [1, 2, 5, 8, 11, 16, 20, 23]). The metric projection operator has been widely applied to nonlinear analysis, optimization theory, approximation theory, fixed point theory, variational analysis, equilibrium theory, and so forth (see [1, 5, 16]).

One of the most important research topics in operator theory is the continuity of the considered operators. Furthermore, in addition to continuity, the smoothness of operators has attracted a lot of attention in the fields of nonlinear analysis. Regarding to the definitions of the smoothness of operators, several types of traditional differentiability of operators, in particular, the metric projection operator, have been introduced and studied in Hilbert spaces (see [3, 5, 6, 9, 10, 13, 14, 19]), and in Banach spaces and normed linear spaces (see [4, 7, 12, 15, 16, 17, 19, 21, 22]).

In this paper, we focus on investigating the Mordukhovich derivatives (see section 2 for the definitions) of the metric projection operator in uniformly convex and uniformly smooth Banach spaces. Therefore, the following concepts of differentiability are only applied to the metric projection operator.

Throughout this paper, unless otherwise stated, let $(X, \|\cdot\|)$ be a real uniformly convex and uniformly smooth Banach space with topological dual space $(X^*, \|\cdot\|_*)$. Let $\langle \cdot, \cdot \rangle$ denote the real canonical pairing between $X^*$ and $X$. Let $C$ be a nonempty closed and convex subset of $X$. Let $P_C: X \to C$ denote the (standard) metric projection operator. We recall some concepts of differentiability of the metric projection operator $P_C$ in uniformly convex and uniformly smooth

Banach spaces below for us to use them in this paper.

I. Gâteaux directional differentiability (see Definition 4.1 in [12]). For $x \in X$ and $v \in X$ with $v \neq \theta$, if the following limit exists,

$$P'_C(x)(v) := \lim_{t \downarrow 0} \frac{P_C(x+tv) - P_C(x)}{t},$$

then, $P_C$ is said to be Gâteaux directionally differentiable at point $x$ along direction $v$. $P'_C(x)(v)$ is called the Gâteaux directional derivative of $P_C$ at $x$ along the direction $v$, and $v$ is called a Gâteaux differentiable direction of $P_C$ at $x$.

II. Fréchet differentiability (see Definition 1.13 in [18]). For any given $\bar{x} \in X$, if there is a linear and continuous mapping $\nabla P_C(\bar{x}): X \to X$ such that

$$\lim_{x \to \bar{x}} \frac{P_C(x) - P_C(\bar{x}) - \nabla P_C(\bar{x})(x - \bar{x})}{\|x - \bar{x}\|} = \theta.$$

then $P_C$ is said to be Fréchet differentiable at $\bar{x}$ and $\nabla P_C(\bar{x})$ is called the Fréchet derivative of $P_C$ at $\bar{x}$.

III. Strict Fréchet differentiability (see Definition 1.13 in [18]). For any given $\bar{x} \in X$, if there is a linear and continuous mapping $\nabla P_C(\bar{x}): X \to X$ such that

$$\lim_{(u,v) \to (\bar{x}, \bar{x})} \frac{P_C(u) - P_C(v) - \nabla P_C(\bar{x})(u-v)}{\|u-v\|} = \lim_{u \to \bar{x}, v \to \bar{x}} \frac{P_C(u) - P_C(v) - \nabla P_C(\bar{x})(u-v)}{\|u-v\|} = \theta,$$

then $P_C$ is said to be strictly Fréchet differentiable at $\bar{x}$ at $\bar{x}$ and $\nabla P_C(\bar{x})$ is called the Fréchet derivative of $P_C$ at $\bar{x}$.

We notice that the above concepts of the differentiability applied to the metric projection operator inherited and perpetuated the definitions of differentiation in functional analysis. Hence, we consider the above concepts of the differentiability of the metric projection operator as the traditional differentiability.

In [18], Mordukhovich successfully introduced some generalized fresh concepts of differentiability of set-valued mappings (Which include single-valued mappings as a special case) in Banach spaces, which has been called generalized differentiation (see Definitions 1.13 and 1.32 in Chapter 1 in [18]) and we call them Mordukhovich differentiation. For the single-valued metric projection operator $P_C$, the coderivative or Fréchet coderivative of $P_C$ at $(\bar{x}, P_C(\bar{x}))$ is denoted by $\widehat{D}^* P_C(\bar{x})$ that is defined to be a set-valued mapping from $X^*$ to $X^*$ such that, for any $y^* \in X^*$,

$$\widehat{D}^* P_C(\bar{x})(y^*) = \left\{ x^* \in X^* : \limsup_{u \to \bar{x}} \frac{\langle x^*, u - \bar{x} \rangle - \langle y^*, P_C(u) - P_C(\bar{x}) \rangle}{\|u - \bar{x}\| + \|P_C(u) - P_C(\bar{x})\|} \leq 0 \right\}.$$

Theorem 1.38 in [18] shows that there are very close connections between the Fréchet derivative $\nabla P_C(\bar{x})$ and the Fréchet coderivative $\widehat{D}^* P_C(\bar{x})$. In particular, if the underlying space is a Hilbert space, then, they will coincide, provided the existence of the Fréchet derivative $\nabla P_C(\bar{x})$. Hence, under the results of Theorem 1.38 in [18], the concept of the Fréchet coderivative $\widehat{D}^* P_C(\bar{x})$ is an extension of the Fréchet derivative $\nabla P_C(\bar{x})$.

In contrast to the traditional differentiability theory (which includes Gâteaux directional differentiability and Fréchet differentiability), the concepts of the (basic) coderivative theory in [18] are new and groundbreaking in nonlinear analysis, which opens a new door in analysis theory and its applications. The coderivative construction introduced in [18] has become a fundamental method of modern analysis theory. It has widely influenced several branches of mathematics such as operator theory, optimization theory, approximation theory, control theory, equilibrium theory, and so forth.

Hence, for clarification (and distinction) of the differences between the concepts of the Fréchet derivative and the Fréchet coderivative, and for paying attention to the importance of the coderivative theory developed in [18], in this paper and in [14], we name the coderivative $\widehat{D}^* P_C(\bar{x})$ as Mordukhovich derivative of the metric projection operator $P_C$ at point $\bar{x}$.

Very recently, in [13, 15], the present author studied the strict Fréchet differentiability of the metric projection operator in Hilbert spaces and the Fréchet differentiability of the metric projection operator in uniformly convex and uniformly smooth Banach spaces, respectively. In [14], the present author investigated Mordukhovich derivatives, in which the considered operator is the metric projection operator in Hilbert spaces. In this paper, we extend the results obtained in [14] to uniformly convex and uniformly smooth Banach spaces. When the considered subset $C$ is a ball, or a closed and convex cylinder in real Banach space $l_p$, or the positive cone in real Banach space $L_p$, we find some precise Mordukhovich derivatives of the metric projection operator.

## 2. Preliminaries

### 2.1. The normalized duality mapping in Banach spaces

Let $(X, \|\cdot\|)$ be a real uniformly convex and uniformly smooth Banach space with topological dual space $(X^*, \|\cdot\|_*)$. Let $\langle \cdot, \cdot \rangle$ denote the real canonical pairing between $X^*$ and $X$. Let $\theta$ and $\theta^*$ denote the origins in $X$ and $X^*$, respectively. Let $\mathbb{B}$ and $\mathbb{B}^*$ denote the unit closed balls in $X$ and $X^*$, respectively. It follows that, for any $r > 0$, $r\mathbb{B}$ and $r\mathbb{B}^*$ are closed balls with radius $r$ and centered at the origins in $X$ and $X^*$, respectively. Let $\mathbb{S}$ be the unit sphere in $X$. Then, $r\mathbb{S}$ is the sphere in $X$ with radius $r$ and center $\theta$. For any $c \in X$ and $r > 0$, let $\mathbb{B}(c, r)$ denote the closed ball in $X$ with radius $r$ and center $c$. It follows that $\mathbb{B}(\theta, 1) = \mathbb{B}$ and $\mathbb{B}(\theta, r) = r\mathbb{B}$. The identity mappings on $X$ and $X^*$ are respectively denoted by $I_X$ and $I_{X^*}$. Let $J: X \to X^*$ and $J^*: X^* \to X$ be the normalized duality mappings. We have

(i) $\langle J(x), x \rangle = \|x\| \|J(x)\|_* = \|x\|^2 = \|J(x)\|_*^2$, for any $x \in X$;

(ii) $\langle x^*, J^*(x^*) \rangle = \|J^*(x^*)\| \|x^*\|_* = \|J^*(x^*)\|^2 = \|x^*\|_*^2$, for any $x^* \in X^*$.

The normalized duality mapping $J$ of the considered uniformly convex and uniformly smooth Banach space $X$ has the following useful properties.

**Lemma 2.1 in [15]**. *Let X be a uniformly convex and uniformly smooth Banach space. For any*

$x, y \in X$, one has

$$2\langle J(y), x - y\rangle \leq \|x\|^2 - \|y\|^2 \leq 2\langle J(x), x - y\rangle.$$

The norm of $X$ is Gâteaux directionally differentiable at every point (see [1, 2, 15, 23]). That is, the following limit exists

$$\lim_{t\downarrow 0}\frac{\|x+ty\| - \|x\|}{t}, \quad \text{uniformly for } (x, y) \in \mathbb{S} \times \mathbb{S}.$$

Then, for the considered uniformly convex and uniformly smooth Banach space $X$, we define the function of smoothness $\Psi$ of $X$ by

$$\Psi(x, y) = \lim_{t\downarrow 0}\frac{\|x+ty\| - \|x\|}{t}, \quad \text{for any } (x, y) \in X \times X.$$

When $(x, y) \in X \times X$, $\Psi(x, y)$ is written as $\psi(x, y)$.

**Lemma 2.2 in [15].** *For any $x, y \in X$ with $x \neq \theta$, one has*

$$\Psi(x, y) = \lim_{t\downarrow 0}\frac{\|x+ty\| - \|x\|}{t} = \frac{\langle J(x), y\rangle}{\|x\|}. \tag{2.1}$$

*In particular, for $(x, y) \in \mathbb{S} \times \mathbb{S}$, by (2.1), one has*

$$\psi(x, y) = \lim_{t\downarrow 0}\frac{\|x+ty\| - \|x\|}{t} = \langle J(x), y\rangle, \quad \text{uniformly for } (x, y) \in \mathbb{S} \times \mathbb{S}. \tag{2.2}$$

### 2.2. The metric projection in uniformly convex and uniformly smooth Banach spaces

Let $C$ be a nonempty closed and convex subset of this uniformly convex and uniformly smooth Banach space $X$. Let $P_C: X \to C$ denote the (standard) metric projection operator. For any $x \in X$, $P_C x \in C$ such that

$$\|x - P_C x\| \leq \|x - z\|, \text{ for all } z \in C.$$

The metric projection operator $P_C$ has the following useful properties (see [1–2, 11, 23] for more details).

**Proposition 2.6 in [2]** *Let $X$ be a uniformly convex and uniformly smooth Banach space and $C$ a nonempty closed and convex subset of $X$. Then the metric projection $P_C: X \to C$ holds the following properties.*

(i) *The operator $P_C$ is fixed on $C$; that is, $P_C(x) = x$, for any $x \in C$;*

(ii) *$P_C$ has the basic variational properties. For any $x \in X$ and $u \in C$, we have*

$$u = P_C(x) \iff \langle J(x-u), u-z\rangle \geq 0, \text{ for all } z \in C; \tag{2.3}$$

(iii) *$P_C$ is uniformly continuous on each bounded subset in $X$.*

In particular, if $M$ is a closed subspace of $X$, then, the basic variational principle of the metric projection operator $P_M: X \to M$ becomes to the following equality version.

**Lemma 2.3 in [15]**. *Let $M$ be a closed subspace of $X$. For any $x \in X$ and $u \in M$,*

$$u = P_M x \quad \Leftrightarrow \quad \langle J(x-u), u-z \rangle = 0, \text{ for all } z \in M,$$

$$\Leftrightarrow \quad \langle J(x-u), z \rangle = 0, \text{ for all } z \in M,$$

$$\Leftrightarrow \quad J(x-u) \in M^\perp. \tag{2.4}$$

**Lemma 2.4 in [15]**. *Let $M$ be a closed subspace of $X$. For arbitrary $x \in X$ and $y \in M$, we have*

$$P_M(\alpha x + \beta y) = \alpha P_M(x) + \beta y, \text{ for any } \alpha, \beta \in \mathbb{R}. \tag{2.5}$$

*In particular, by taking $y = \theta$ in (2.5), we have*

$$P_M(\alpha x) = \alpha P_M(x), \text{ for any } \alpha \in \mathbb{R}. \tag{2.6}$$

### 2.3. Semi-orthogonal decompositions in Banach spaces

In [15], we introduced the semi-orthogonal decompositions of points in Banach spaces with respect to a given nonzero point, which is repeatedly applied in the proofs of Theorem 3.1 in [15]. In this paper, for studying the Mordukhovich derivatives of the metric projection operator in uniformly convex and uniformly smooth Banach spaces, we extend the semi-orthogonal decompositions of points in Banach spaces to their dual spaces. We first recall the semi-orthogonal decompositions of points in Banach spaces introduced in [15]. For any given $\bar{x} \in X \setminus \{\theta\}$, every $x \in X$ can be clearly written as

$$x = \frac{\langle J(\bar{x}), x \rangle}{\|\bar{x}\|^2} \bar{x} + \left( x - \frac{\langle J(\bar{x}), x \rangle}{\|\bar{x}\|^2} \bar{x} \right), \text{ for all } x \in X. \tag{2.7}$$

This implies

$$\langle J(\bar{x}), x - \frac{\langle J(\bar{x}), x \rangle}{\|\bar{x}\|^2} \bar{x} \rangle = 0, \text{ for all } x \in X.$$

Hence, (2.7) is called a semi-orthogonal decomposition of $x \in X$ with respect to the arbitrarily given $\bar{x} \in X \setminus \{\theta\}$. We define a real valued functional $a(\bar{x}; \cdot): X \to \mathbb{R}$ by

$$a(\bar{x}; x) = \frac{\langle J(\bar{x}), x \rangle}{\|\bar{x}\|^2}, \text{ for all } x \in X.$$

Let

$$O(\bar{x}) := \{ y \in X : \langle J(\bar{x}), y \rangle = 0 \}$$

We define a mapping $o(\bar{x}; \cdot): X \to O(\bar{x})$ by

$$o(\bar{x}; x) := x - \frac{\langle J(\bar{x}), x \rangle}{\|\bar{x}\|^2} \bar{x} = x - a(\bar{x}; x)\bar{x}, \text{ for all } x \in X.$$

By (2.7), with respect to this arbitrarily given $\bar{x} \in X \setminus \{\theta\}$, for every $x \in X$, the semi-orthogonal decomposition (2.7) can be rewritten as

$$x = a(\bar{x}; x)\bar{x} + o(\bar{x}; x), \text{ for all } x \in X. \tag{2.7}$$

The following lemma provides some properties of $a(\bar{x}; \cdot)$ and $o(\bar{x}; \cdot)$. These properties will play important roles and will be repeatedly used in the following section of this paper.

**Proposition 2.11 (partial) in [15].** *For any fixed $\bar{x} \in X\backslash\{\theta\}$, the real valued function $a(\bar{x}; \cdot)$ and the mapping $o(\bar{x}; \cdot)$ have the following properties.*

(i)    $a(\bar{x}; \cdot): X \to \mathbb{R}$ *is a real valued linear and continuous functional satisfying*

$$\langle J(\bar{x}), x \rangle = a(\bar{x}; x)\|\bar{x}\|^2, \quad \text{for any } x \in X;$$

(ii)    $o(\bar{x}; \cdot): X \to O(\bar{x})$ *is a linear and continuous mapping satisfying*

$$\langle J(\bar{x}), o(\bar{x}; x) \rangle = 0, \quad \text{for any } x \in X;$$

(iii)    $u \to \bar{x} \iff a(\bar{x}; u) \to 1 \text{ and } o(\bar{x}; u) \to \theta, \text{ for } u \in X.$

For any given $\bar{x} \in X\backslash\{\theta\}$, we have $J(\bar{x}) \in X^*\backslash\{\theta^*\}$. Let $\bar{x}^* := J(\bar{x}) \in X^*\backslash\{\theta^*\}$. Similar to the semi-orthogonal decomposition (2.7) for points in $X$ with respect to $\bar{x} \in X\backslash\{\theta\}$, every $x^* \in X^*$ has the following semi-orthogonal decomposition with respect to $\bar{x}^* = J(\bar{x}) \in X^*\backslash\{\theta^*\}$.

$$x^* = \frac{\langle x^*, J^*\bar{x}^* \rangle}{\|\bar{x}^*\|_*^2}\bar{x}^* + \left(x^* - \frac{\langle x^*, J^*\bar{x}^* \rangle}{\|\bar{x}^*\|_*^2}\bar{x}^*\right)$$

$$= \frac{\langle x^*, \bar{x} \rangle}{\|\bar{x}^*\|_*^2}\bar{x}^* + \left(x^* - \frac{\langle x^*, \bar{x} \rangle}{\|\bar{x}^*\|_*^2}\bar{x}^*\right)$$

$$= \frac{\langle x^*, \bar{x} \rangle}{\|\bar{x}\|^2}\bar{x}^* + \left(x^* - \frac{\langle x^*, \bar{x} \rangle}{\|\bar{x}\|^2}\bar{x}^*\right), \text{ for all } x^* \in X^*. \tag{2.8}$$

By $\bar{x}^* := J(\bar{x})$ and $\langle \bar{x}^*, \bar{x} \rangle = \|\bar{x}^*\|_*^2 = \|\bar{x}\|^2$, one easily checks that

$$\langle x^* - \frac{\langle x^*, \bar{x} \rangle}{\|\bar{x}\|^2}\bar{x}^*, \bar{x} \rangle = 0, \text{ for all } x^* \in X^*.$$

Hence, (2.8) is called the semi-orthogonal decomposition of $x^* \in X^*$ with respect to $\bar{x}^* = J(\bar{x})$, which is induced by the arbitrarily given $\bar{x} \in X\backslash\{\theta\}$. We similarly define a real valued functional $a^*(\bar{x}^*; \cdot): X^* \to \mathbb{R}$ by

$$a^*(\bar{x}^*; x^*) = \frac{\langle x^*, \bar{x} \rangle}{\|\bar{x}\|^2} = \frac{\langle x^*, \bar{x} \rangle}{\|\bar{x}^*\|_*^2}, \text{ for all } x^* \in X^*.$$

We define a mapping $o^*(\bar{x}^*; \cdot): X^* \to X^*$ by

$$o^*(\bar{x}^*; x^*) := x^* - \frac{\langle x^*, \bar{x} \rangle}{\|\bar{x}\|^2}\bar{x}^* = x^* - a^*(\bar{x}^*; x^*)\bar{x}^*, \text{ for all } x^* \in X^*.$$

Then, the semi-orthogonal decomposition (2.8) of $x^* \in X^*$ with respect to $\bar{x}^* = J(\bar{x})$ induced by the given $\bar{x} \in X\backslash\{\theta\}$ has the following version of the semi-orthogonal decomposition

$$x^* = a^*(\bar{x}^*; x^*)\bar{x}^* + o^*(\bar{x}^*; x^*), \text{ for all } x^* \in X^*. \tag{2.8}$$

Similar to Proposition 2.11 in [15], the following proposition provides properties of $a^*(\bar{x}^*; \cdot)$ and $o^*(\bar{x}^*; \cdot)$. All properties proved in Proposition 2.11 in [15] and the following properties will play important roles in the following section of this paper.

**Proposition 2.1.** *For any fixed $\bar{x} \in X\setminus\{\theta\}$, let $\bar{x}^* = J(\bar{x}) \in X^*\setminus\{\theta^*\}$. The real valued function $a^*(\bar{x}^*; \cdot)$ and the mapping $o^*(\bar{x}^*; \cdot)$ have the following properties.*

(i) $a^*(\bar{x}^*; \cdot): X^* \to \mathbb{R}$ *is a real valued linear and continuous functional satisfying*

$$\langle x^*, \bar{x}\rangle = a^*(\bar{x}^*; x^*)\|\bar{x}\|^2 = a^*(\bar{x}^*; x^*)\|\bar{x}^*\|_*^2, \text{ for any } x^* \in X^*;$$

(ii) $o^*(\bar{x}^*; \cdot): X^* \to X^*$ *is a linear and continuous mapping satisfying*

$$\langle o^*(\bar{x}^*; x^*), \bar{x}\rangle = 0, \text{ for any } x^* \in X^*;$$

(iii) $u^* \to \bar{x}^* \iff a^*(\bar{x}^*; u^*) \to 1 \text{ and } o^*(\bar{x}^*; u^*) \to \theta^*, \text{ for } u^* \in X^*;$

All $a(\bar{x}; \cdot)$, $o(\bar{x}; \cdot)$, $a^*(\bar{x}^*; \cdot)$ and $o^*(\bar{x}^*; \cdot)$ depend (directly or indirectly) on $\bar{x}$. However, for the sake of simplicity, $a(\bar{x}; \cdot)$, $o(\bar{x}; \cdot)$, $a^*(\bar{x}^*; \cdot)$ and $o^*(\bar{x}^*; \cdot)$ are abbreviated as $a(\cdot)$, $o(\cdot)$, $a^*(\cdot)$ and $o^*(\cdot)$, respectively.

**Proposition 2.12 in [15].** *Let $(X, \|\cdot\|)$ be a real uniformly convex and uniformly smooth Banach space. Then, we have*

$$\lim_{\substack{v \to \theta \\ v \in O(\bar{x})}} \frac{\|\bar{x}+v\| - \|\bar{x}\|}{\|v\|} = 0, \text{ for each } \bar{x} \in X\setminus\{\theta\}. \tag{2.10}$$

### 2.4. The product spaces of Banach spaces

Let $X \times X$ be the Cartesian product space of $X$. In this paper, we define norms on $X \times X$ and $X^* \times X^*$, which are denoted by $\|\cdot\|_A$ and $\|\cdot\|_{A^*}$, respectively. We also define the real canonical pairing between $X^* \times X^*$ and $X \times X$. So that the normalized duality mapping on $X \times X$ will be defined with respect to the norms $\|\cdot\|_A$ and $\|\cdot\|_{A^*}$ and the canonical pairing.

In this paper, the norm $\|\cdot\|_A$ on $X \times X$ is defined by

$$\|(x, y)\|_A = \sqrt{\|x\|^2 + \|y\|^2}, \text{ for all } (x, y) \in X \times X. \tag{2.11}$$

Similarly, the norm $\|\cdot\|_{A^*}$ on $X^* \times X^*$ is defined by

$$\|(x^*, y^*)\|_{A^*} = \sqrt{\|x^*\|_*^2 + \|x^*\|_*^2}, \text{ for all } (x^*, y^*) \in X^* \times X^*. \tag{2.12}$$

Let $\langle \cdot, \cdot \rangle_A$ denote the real canonical pairing between $X^* \times X^*$ and $X \times X$, which is defined by

$$\langle (x^*, y^*), (x, y)\rangle_A = \langle x^*, x\rangle + \langle y^*, y\rangle, \text{ for all } (x, y) \in X \times X \text{ and } (x^*, y^*) \in X^* \times X^*. \tag{2.13}$$

Based on $J: X \to X^*$, we define a mapping $\mathbb{J}: X \times X \to X^* \times X^*$ by

$$\mathbb{J}(x,y) := (Jx, Jy) \in X^* \times X^*, \text{ for any } (x,y) \in X \times X. \tag{2.14}$$

By (2.13) and (2.11), for all $(x,y) \in X \times X$, we have

$$\langle \mathbb{J}(x,y), (x,y) \rangle_A$$
$$= \langle (Jx, Jy), (x,y) \rangle_A$$
$$= \langle Jx, x \rangle + \langle Jy, y \rangle$$
$$= \|x\|^2 + \|y\|^2$$
$$= \|(x,y)\|_A^2. \tag{2.15}$$

On the other hand, by (2.13) and (2.12) and by the definition of $J$, we have

$$\langle \mathbb{J}(x,y), (x,y) \rangle_A$$
$$= \|x\|^2 + \|y\|^2$$
$$= \|Jx\|_*^2 + \|Jy\|_*^2$$
$$= \|(Jx, Jy)\|_{A^*}^2. \tag{2.16}$$

By (2.15) and (2.16), we obtain that $\mathbb{J}$ is the normalized duality mapping from $X \times X$ to $X^* \times X^*$. The Cartesian product space $X \times X$ with the norm $\|\cdot\|_A$ is also a uniformly convex and uniformly smooth Banach space with the dual space $X^* \times X^*$.

As what is mentioned in the proof of Proposition 12 in Chapter 1 of [18], the prenormal cone $\widehat{N}(x; \Omega)$ (which will be recalled in the following subsection) does not depend on equivalent norms on $X$. In [18], the norm on $X \times X$ is defined by

$$\|(x,y)\| = \|x\| + \|y\|, \text{ for all } (x,y) \in X \times X. \tag{2.17}$$

Notice that

$$\tfrac{\sqrt{2}}{2}(\|x\| + \|y\|) \leq \sqrt{\|x\|^2 + \|y\|^2} \leq \|x\| + \|y\|, \text{ for all } (x,y) \in X \times X. \tag{2.18}$$

This implies that the norm $\|\cdot\|_A$ defined in (2.11) is equivalent to the norm defined in (2.17) on the Cartesian product space $X \times X$. Hence, in the following sections, we use the norm on the product space $X \times X$ defined in (2.11).

### 2.5. Mordukhovich derivatives of the metric projection operator in Banach spaces

In [18], the concepts of Mordukhovich derivatives of set-valued mappings in Banach spaces are introduced. In this paper, we focus on investigating Mordukhovich derivatives of the metric projection operator in uniformly convex and uniformly smooth Banach spaces. Here, we rewrite the definitions of Mordukhovich derivatives given in [18] with respect to this single-valued mapping $P_C$. The definition of generalized differentiation is based on the prenormal cones.

**Definition 1.1 in [18]** (generalized normals in uniformly convex and uniformly smooth Banach spaces). Let $\Omega$ be a nonempty subset of a uniformly convex and uniformly smooth Banach space $X$. For any $x \in X$, we define

$$\widehat{N}(x;\Omega) = \left\{x^* \in X^* : \limsup_{\substack{u \to x \\ \Omega}} \frac{\langle x^*, u-x \rangle}{\|u-x\|} \leq 0\right\}, \text{ for any } x \in \Omega. \tag{2.19}$$

The elements of (2.9) are called Fréchet normals and $\widehat{N}(x;\Omega)$ is the prenormal cone to $\Omega$ at $x$. We put $\widehat{N}(x;\Omega) = \emptyset$, for any $x \in X\backslash\Omega$. Then, based on the above definition, we recall the Mordukhovich derivatives of the single-valued mapping $P_C$ in uniformly convex and uniformly smooth Banach spaces. Note that, in [18], Mordukhovich derivatives are named by Fréchet coderivatives or precoderivatives.

**Definition 1.32 in [18]** (Mordukhovich derivatives for single-valued mapping $P_C$ in uniformly convex and uniformly smooth Banach spaces).

The Mordukhovich derivatives of $P_C$ at $(x, P_C(x))$ is defined, for any given $y^* \in X^*$, by

$$\widehat{D}^* P_C(x, P_C(x))(y^*) := \widehat{D}^* P_C(x)(y^*) = \{x^* \in X^* : (x^*, -y^*) \in \widehat{N}((x, P_C(x)); \mathrm{gph} P_C)\}.$$

It follows that, if $u \neq P_C(x)$, then

$$\widehat{D}^* P_C(x, u)(y^*) = \emptyset, \text{ for all } y^* \in X^*.$$

Next theorem provides the connections between Mordukhovich derivatives and Fréchet derivatives for single-valued mappings in uniformly convex and uniformly smooth Banach spaces.

**Theorem 1.38 in [18].** *Suppose that $P_C: X \to C$ is Fréchet differentiable at $\bar{x} \in X$. Then, Mordukhovich derivative of $P_C$ at $\bar{x}$ satisfies the following equation*

$$\widehat{D}^* P_C(\bar{x})(y^*) = \{(\nabla P_C(\bar{x}))^*(y^*)\}, \text{ for all } y^* \in X^*.$$

By the continuity of $P_C$ in in uniformly convex and uniformly smooth Banach spaces, we calculate the solutions of the Mordukhovich derivative of $P_C$ at $(x, P_C(x))$ as follows. For any $y^* \in X^*$, we have

$$\widehat{D}^* P_C(x)(y^*)$$
$$= \{x^* \in X^* : (x^*, -y^*) \in \widehat{N}((x, P_C(x)); \mathrm{gph} P_C)\}$$
$$= \left\{x^* \in X^* : \limsup_{(u, P_C(u)) \to (x, P_C(x))} \frac{\langle (x^*, -y^*), (u, P_C(u)) - (x, P_C(x)) \rangle_A}{\|(u, P_C(u)) - (x, P_C(x))\|_A} \leq 0\right\}$$
$$= \left\{x^* \in X^* : \limsup_{u \to x} \frac{\langle x^*, u-x \rangle - \langle y^*, P_C(u) - P_C(x) \rangle}{\sqrt{\|u-x\|^2 + \|P_C(u) - P_C(x)\|^2}} \leq 0\right\}$$

$$= \left\{ x^* \in X^* : \limsup_{u \to x} \frac{\langle x^*, u-x \rangle - \langle y^*, P_C(u)-P_C(x) \rangle}{\|u-x\| + \|P_C(u)-P_C(x)\|} \leq 0 \right\}. \qquad (2.20)$$

## 3. Mordukhovich derivatives of the metric projection onto closed balls

Theorem 3.3 in [13] proves the strict Fréchet differentiability of the metric projection onto closed balls in Hilbert spaces, which is applied to investigate Mordukhovich derivatives of the metric projection in Hilbert spaces. Theorem 3.1 in [15] proves the Fréchet differentiability of the metric projection onto closed balls in uniformly convex and uniformly smooth Banach spaces. However, in the following theorem, we investigate Mordukhovich derivatives of the metric projection onto closed balls in uniformly convex and uniformly smooth Banach spaces and find its exact representations of the Mordukhovich derivatives. Before we prove next theorem, we review some properties of the metric projection onto closed balls. For any $r > 0$, the metric projection $P_{r\mathbb{B}}: X \to r\mathbb{B}$ has the following analytic representations (see [15]).

$$P_{r\mathbb{B}}(x) = \begin{cases} x, & \text{for any } x \in r\mathbb{B}, \\ \frac{r}{\|x\|} x, & \text{for any } x \notin r\mathbb{B}. \end{cases}$$

We need the following notations for the following theorem. For any $x \in r\mathbb{S}$, two subsets $x_r^\uparrow$ and $x_r^\downarrow$ of $X\setminus\{\theta\}$ are defined by

(a) $x_r^\uparrow = \{v \in X\setminus\{\theta\}: \text{there is } \delta > 0 \text{ such that } \|x + tv\| > r, \text{ for all } t \in (0, \delta)\}$;
(b) $x_r^\downarrow = \{v \in X\setminus\{\theta\}: \text{there is } \delta > 0 \text{ such that } \|x + tv\| \leq r, \text{ for all } t \in (0, \delta)\}$.

**Theorem 3.1 (partial) in [15].** *Let $X$ be a uniformly convex and uniformly smooth Banach space. For any $r > 0$, the metric projection $P_{r\mathbb{B}}: X \to r\mathbb{B}$ has the following differentiable properties.*

(i) $P_{r\mathbb{B}}$ *is strictly Fréchet differentiable on $r\mathbb{B}^o$ satisfying*

$$\nabla P_{r\mathbb{B}}(\bar{x}) = I_X, \text{ for every } \bar{x} \in r\mathbb{B}^o.$$

*That is,*

$$\bar{x} \in r\mathbb{B}^o \implies \nabla P_{r\mathbb{B}}(\bar{x})(x) = x, \text{ for every } x \in X;$$

(ii) $P_{r\mathbb{B}}$ *is Fréchet differentiable at every $\bar{x} \in X\setminus r\mathbb{B}$. The Fréchet derivative at $\bar{x}$ satisfies*

$$\nabla P_{r\mathbb{B}}(\bar{x})(x) = \frac{r}{\|\bar{x}\|}\left(x - \frac{\langle J(\bar{x}), x \rangle}{\|\bar{x}\|^2} \bar{x}\right), \text{ for every } x \in X.$$

In next theorem, we use Theorem 3.1 in [15] to investigate the solutions and properties of Mordukhovich derivatives of the metric projection $P_{r\mathbb{B}}: X \to r\mathbb{B}$, which is an extension of Theorem 3.2 in [14] from Hilbert spaces to uniformly convex and uniformly smooth Banach spaces. The results and the proof of the following theorem are similar to that in Theorem 3.2 in [14]. However, in the proof of the following theorem, we use the semi-decompositions of Banach spaces and their dual spaces proved in Proposition 2.1.

**Theorem 3.1.** *Let X be a uniformly convex and uniformly smooth Banach space. For any $r > 0$, the Mordukhovich derivatives of the metric projection $P_{r\mathbb{B}}\colon X \to r\mathbb{B}$ has the following properties.*

(i)   *For every $\bar{x} \in r\mathbb{B}^o$, we have*
$$\widehat{D}^*P_C(\bar{x})(y^*) = \{y^*\}, \text{ for every } y^* \in X^*.$$

(ii)  *For every $\bar{x} \in X \setminus r\mathbb{B}$, we have*
$$\widehat{D}^*P_C(\bar{x})(y^*) = \left\{\frac{r}{\|\bar{x}\|}\left(y^* - \frac{\langle y^*, x\rangle}{\|\bar{x}\|^2}J(\bar{x})\right)\right\}, \text{ for every } y^* \in X^*.$$

*In particular, we have*

(a) $\widehat{D}^*P_C(\bar{x})(y^*) = \left\{\frac{r}{\|\bar{x}\|}y^*\right\}$, *if* $y^* \perp \bar{x}$;

(b) $\widehat{D}^*P_C(\bar{x})(J(\bar{x})) = \{\theta^*\}$.

(iii) *If $\bar{x} \in r\mathbb{S}$, then*

(a)   $\widehat{D}^*P_{r\mathbb{B}}(\bar{x})(\theta^*) = \{\theta^*\}$;

(b)   *For any $y^* \in X^* \setminus \{\theta^*\}$, we have*

   *Part 1.*
   $$\theta^* \in \widehat{D}^*P_{r\mathbb{B}}(\bar{x})(y^*)$$
   $$\Rightarrow \langle y^*, \bar{x}\rangle \leq 0, -J^*(o^*(y^*)\bar{x}^*) \notin \bar{x}_r^{\uparrow}$$
   $$\text{and } \frac{\langle y^*, \bar{x}\rangle}{r^2}\langle \bar{x}^*, J^*(o^*(y^*))\rangle + \|J^*(o^*(y^*))\|^2 \leq 0.$$

   *Part 2.*
   $$o^*(y^*) = \theta^* \text{ and } \langle y^*, \bar{x}\rangle < 0 \Rightarrow \theta^* \in \widehat{D}^*P_{r\mathbb{B}}(\bar{x})(y^*).$$

   *In particular,* $\theta^* \in \widehat{D}^*P_{r\mathbb{B}}(\bar{x})(-J(\bar{x}))$.

(c)   $\widehat{D}^*P_{r\mathbb{B}}(\bar{x})(J(\bar{x})) = \emptyset.$

Before we prove this theorem, we note that, comparing the special case in part (b) in (iii) to part (c) in (iii), we find that, for any given $\bar{x} \in r\mathbb{S}$, $\widehat{D}^*P_{r\mathbb{B}}(\bar{x})\colon X^* \to X^*$ is not a linear mapping, in general.

*Proof.* Proof of (i). By part (i) in Theorem 3.1 in [14], we have
$$\nabla P_{r\mathbb{B}}(\bar{x}) = I_X, \text{ for every } \bar{x} \in r\mathbb{B}^o.$$

Then, by Theorem 1.38 in [18], for every $\bar{x} \in r\mathbb{B}^o$, we have

$$\widehat{D}^*P_C(\bar{x})(y^*) = \{(I_X)^*(y^*)\} = \{y^*\}, \text{ for all } y^* \in X^*.$$

Proof of (ii). By part (ii) in Theorem 3.1 in [14], for every $\bar{x} \in X \backslash r\mathbb{B}$, we have

$$\nabla P_{r\mathbb{B}}(\bar{x})(x) = \frac{r}{\|\bar{x}\|}\left(x - \frac{\langle J(\bar{x}),x\rangle}{\|\bar{x}\|^2}\bar{x}\right), \text{ for every } x \in X.$$

Then, by Theorem 1.38 in [18], for every $\bar{x} \in X \backslash r\mathbb{B}$, we have

$$\widehat{D}^*P_C(\bar{x})(y^*) = \{(\nabla P_C(\bar{x}))^*(y^*)\}, \text{ for all } y^* \in X^*.$$

Since $(\nabla P_C(\bar{x}))^*(y^*) \in X^*$, for any $x \in X$, we have

$$\langle (\nabla P_C(\bar{x}))^*(y^*), x \rangle$$

$$= \langle y^*, \nabla P_{r\mathbb{B}}(\bar{x})(x) \rangle$$

$$= \langle y^*, \frac{r}{\|\bar{x}\|}\left(x - \frac{\langle J(\bar{x}),x\rangle}{\|\bar{x}\|^2}\bar{x}\right) \rangle$$

$$= \frac{r}{\|\bar{x}\|}\langle y^*, x - \frac{\langle J(\bar{x}),x\rangle}{\|\bar{x}\|^2}\bar{x} \rangle$$

$$= \frac{r}{\|\bar{x}\|}\left(\langle y^*, x\rangle - \langle y^*, \frac{\langle J(\bar{x}),x\rangle}{\|\bar{x}\|^2}\bar{x}\rangle\right)$$

$$= \frac{r}{\|\bar{x}\|}\left(\langle y^*, x\rangle - \frac{\langle y^*,\bar{x}\rangle}{\|\bar{x}\|^2}\langle J(\bar{x}), x\rangle\right)$$

$$= \langle \frac{r}{\|\bar{x}\|}\left(y^* - \frac{\langle y^*,\bar{x}\rangle}{\|\bar{x}\|^2}J(\bar{x})\right), x\rangle, \text{ for any } x \in X.$$

This implies

$$(\nabla P_C(\bar{x}))^*(y^*) = \frac{r}{\|\bar{x}\|}\left(y^* - \frac{\langle y^*,\bar{x}\rangle}{\|\bar{x}\|^2}J(\bar{x})\right).$$

This proves (ii).

Proof of (a) in (iii). Since $\|\bar{x}\| = r$, we have $P_{r\mathbb{B}}(\bar{x}) = \bar{x}$. It is clear to see

$$\theta^* \in \widehat{D}^*P_{r\mathbb{B}}(\bar{x})(\theta^*). \tag{3.1}$$

For any $z^* \in X^*\backslash\{\theta^*\}$, we want to show $z^* \notin \widehat{D}^*P_{r\mathbb{B}}(\bar{x})(\theta^*)$. By (2.20), we have

$$z^* \in \widehat{D}^*P_{r\mathbb{B}}(\bar{x})(\theta^*) \iff \limsup_{u \to \bar{x}} \frac{\langle(z^*,-\theta^*),\ (u,P_{r\mathbb{B}}(u)-(\bar{x},P_{r\mathbb{B}}(\bar{x})))\rangle}{\|u-\bar{x}\|+\|P_{r\mathbb{B}}(u)-P_{r\mathbb{B}}(\bar{x})\|} \leq 0$$

$$\iff \limsup_{u \to \bar{x}} \frac{\langle(z^*,-\theta^*),\ (u,P_{r\mathbb{B}}(u)-(\bar{x},\bar{x}))\rangle}{\|u-\bar{x}\|+\|P_{r\mathbb{B}}(u)-\bar{x}\|} \leq 0$$

$$\Leftrightarrow \quad \limsup_{u \to \bar{x}} \frac{\langle z^*, u-\bar{x}\rangle}{\|u-\bar{x}\|+\|P_{r\mathbb{B}}(u)-\bar{x}\|} \leq 0. \tag{3.2}$$

Let $z = J^*(z^*)$. Since $z^* \in X^*\setminus\{\theta^*\}$, then

$$\langle z^*, z\rangle = \|z\|^2 = \|z^*\|_*^2 > 0.$$

Case 1. Suppose $z \in \bar{x}_r^{\uparrow}$. Then, there is $\delta > 0$ such that $\|\bar{x}+tz\| > r$, for all $t \in (0, \delta)$. In this case, we take a directional line segment in the limit (3.2) as $u = \bar{x} + tz$, for $t \downarrow 0$ with $t \in (0, \delta)$. We have

$$z^* \in \widehat{D}^* P_{r\mathbb{B}}(\bar{x})(\theta^*) \quad \Leftrightarrow \quad 0 \geq \limsup_{u \to \bar{x}} \frac{\langle z^*, u-\bar{x}\rangle}{\|u-\bar{x}\|+\|P_{r\mathbb{B}}(u)-\bar{x}\|}$$

$$\geq \limsup_{t \downarrow 0, t < \delta} \frac{\langle z^*, \bar{x}+tz-\bar{x}\rangle}{\|\bar{x}+tz-\bar{x}\|+\|P_{r\mathbb{B}}(\bar{x}+tz)-\bar{x}\|}$$

$$= \limsup_{t \downarrow 0, t < \delta} \frac{t\langle z^*, z\rangle}{t\|z\|+\|P_{r\mathbb{B}}(\bar{x}+tz)-\bar{x}\|}$$

$$= \limsup_{t \downarrow 0, t < \delta} \frac{t\|z\|^2}{t\|z\|+\left\|\frac{r}{\|\bar{x}+tz\|}(\bar{x}+tz)-\bar{x}\right\|}. \tag{3.6}$$

We calculate $\left\|\frac{r}{\|\bar{x}+tz\|}(\bar{x}+tz) - \bar{x}\right\|$ with $t \in (0, \delta)$ in the denominator in limit (3.6).

$$\left\|\frac{r}{\|\bar{x}+tz\|}(\bar{x}+tz) - \bar{x}\right\|$$

$$\leq \|\bar{x}\|\left|\frac{r}{\|\bar{x}+tz\|}-1\right| + \frac{tr}{\|\bar{x}+tz\|}\|z\|$$

$$= r\left|\frac{\|\bar{x}+tz\|-r}{\|\bar{x}+tz\|}\right| + \frac{tr}{\|\bar{x}+tz\|}\|z\|$$

$$= \frac{tr}{\|\bar{x}+tz\|}\left|\frac{\|\bar{x}+tz\|-r}{t}\right| + \frac{tr\|z\|}{\|\bar{x}+tz\|}$$

$$= \frac{tr}{\|\bar{x}+tz\|}\left|\frac{\|\bar{x}+tz\|-\|\bar{x}\|}{t}\right| + \frac{tr\|z\|}{\|\bar{x}+tz\|}$$

$$= \frac{tr}{\|\bar{x}+tz\|}\left(\left|\frac{\|\bar{x}+tz\|-\|\bar{x}\|}{t}\right| + \|z\|\right). \tag{3.7}$$

Substituting (3.7) into the limit (3.6), by (2.1) and (3.5), we have

$$0 \geq \limsup_{u \to \bar{x}} \frac{\langle z^*, u-\bar{x}\rangle}{\|u-\bar{x}\|+\|P_{r\mathbb{B}}(u)-\bar{x}\|}$$

$$\geq \limsup_{t \downarrow 0, t < \delta} \frac{t\langle z^*, z\rangle}{t\|z\|+\left\|\frac{r}{\|\bar{x}+tz\|}(\bar{x}+tz)-\bar{x}\right\|}$$

$$\geq \limsup_{t\downarrow 0, t<\delta} \frac{t\|z\|^2}{t\|z\|+\frac{tr}{\|\bar{x}+tz\|}\left(\left|\frac{\|\bar{x}+tz\|-\|\bar{x}\|}{t}\right|+\|z\|\right)}$$

$$= \limsup_{t\downarrow 0, t<\delta} \frac{\|z\|^2}{\|z\|+\frac{r}{\|\bar{x}+tz\|}\left(\left|\frac{\|\bar{x}+tz\|-\|\bar{x}\|}{t}\right|+\|z\|\right)}$$

$$= \limsup_{t\downarrow 0, t<\delta} \frac{\|z\|^2}{\|z\|+\frac{r}{\|\bar{x}\|}\left(\left|\frac{\langle J(\bar{x}), z\rangle}{\|\bar{x}\|}\right|+\|z\|\right)}$$

$$= \frac{\|z\|^2}{2\|z\|+\left|\frac{\langle J(\bar{x}), z\rangle}{r}\right|}$$

$$\geq \frac{\|z\|^2}{2\|z\|+\frac{1}{r}\|J(\bar{x})\|_*\|z\|}$$

$$= \frac{\|z\|^2}{2\|z\|+\frac{1}{r}\|\bar{x}\|\|z\|}$$

$$= \frac{\|z\|}{3}$$

$$> 0. \tag{3.8}$$

By (3.6), this contradicts to the statement $z^* \in \widehat{D}^* P_{r\mathbb{B}}(\bar{x})(\theta^*)$.

Case 2. Suppose $z \in \bar{x}_r^{\downarrow}$. Then, there is $\delta > 0$ such that $\|\bar{x}+tz\| \leq r$, for all $t \in (0, \delta)\}$. In this case, we take a directional line segment in the limit (3.2) as, $u = \bar{x} + tz$, for $t \downarrow 0$ with $t \in (0, \delta)$. We have

$$0 \geq \limsup_{u \to \bar{x}} \frac{\langle z^*, u-\bar{x}\rangle}{\|u-\bar{x}\|+\|P_{r\mathbb{B}}(u)-\bar{x}\|}$$

$$= \limsup_{t\downarrow 0, t<\delta} \frac{t\langle z^*, z\rangle}{t\|z\|+\|\bar{x}+tz-\bar{x}\|}$$

$$= \limsup_{t\downarrow 0, t<\delta} \frac{t\|z\|^2}{2t\|z\|}$$

$$= \frac{\|z\|}{2}$$

$$> 0.$$

By (3.6), this contradicts to the statement $z^* \in \widehat{D}^* P_{r\mathbb{B}}(\bar{x})(\theta^*)$. With (3.8), this proves

$$z^* \in X^*\setminus\{\theta^*\} \implies z^* \notin \widehat{D}^* P_{r\mathbb{B}}(\bar{x})(\theta^*).$$

By (3.1) and (3.9), we proved (a) of part (iii), that is,

$$\widehat{D}^* P_{r\mathbb{B}}(\bar{x})(\theta^*) = \{\theta^*\}, \text{ for any } \bar{x} \in r\mathbb{S}.$$

Proof of part 1 of (b) in (iii). Let $y^* \in X^*\setminus\{\theta^*\}$. We show that

$$\theta^* \in \widehat{D}^*P_{r\mathbb{B}}(\bar{x})(y^*)$$

$$\Rightarrow \langle y^*, \bar{x}\rangle \leq 0, -J^*(o^*(y^*)\bar{x}^*) \notin \bar{x}_r^{\uparrow} \text{ and } \frac{\langle y^*, \bar{x}\rangle}{r^2}\langle \bar{x}^*, J^*(o^*(y^*))\rangle + \|J^*(o^*(y^*))\|^2 \leq 0. \quad (3.9)$$

By definition, for $y^* \in X^*$ with $y^* \neq \theta^*$,

$$\theta^* \in \widehat{D}^*P_{r\mathbb{B}}(\bar{x})(y^*) \Leftrightarrow \limsup_{u \to \bar{x}} \frac{\langle (\theta^*,-y^*), (u,P_{r\mathbb{B}}(u)-(\bar{x},P_{r\mathbb{B}}(\bar{x})))\rangle}{\|u-\bar{x}\| + \|P_{r\mathbb{B}}(u)-P_{r\mathbb{B}}(\bar{x})\|} \leq 0$$

$$\Leftrightarrow \limsup_{u \to \bar{x}} \frac{-\langle y^*, P_{r\mathbb{B}}(u)-\bar{x}\rangle}{\|u-\bar{x}\| + \|P_{r\mathbb{B}}(u)-\bar{x}\|} \leq 0. \quad (3.10)$$

We take a directional line segment in the limit in (3.10), $u = (1-t)\bar{x}$, for $t \downarrow 0$ with $t < 1$. Then, by the semi-decomposition (2.8) of the point $y^* \in X^*$ with respect to $\bar{x}^* = J(\bar{x}) \in X^*\setminus\{\theta^*\}$ induced by the given point $\bar{x} \in r\mathbb{S}$, and by Proposition 2.1, we have

$$0 \geq \limsup_{u \to \bar{x}} \frac{-\langle y^*, P_{r\mathbb{B}}(u)-\bar{x}\rangle}{\|u-\bar{x}\|+\|P_{r\mathbb{B}}(u)-\bar{x}\|}$$

$$\geq \limsup_{t\downarrow 0, t<1} \frac{-\langle y^*, P_{r\mathbb{B}}((1-t)\bar{x})-\bar{x}\rangle}{\|(1-t)\bar{x}-\bar{x}\|+\|P_{r\mathbb{B}}((1-t)\bar{x})-\bar{x}\|}$$

$$= \limsup_{t\downarrow 0, t<1} \frac{-\langle a^*(y^*)\bar{x}^*+o^*(y^*), -t\bar{x}\rangle}{\|(1-t)\bar{x}-\bar{x}\|+\|(1-t)\bar{x}-\bar{x}\|}$$

$$= \limsup_{t\downarrow 0, t<1} \frac{t(\langle a^*(y^*)\bar{x}^*, \bar{x}\rangle + \langle o^*(y^*), \bar{x}\rangle)}{2t\|\bar{x}\|}$$

$$= \limsup_{t\downarrow 0, t<1} \frac{ta^*(y^*)\|\bar{x}\|^2}{2t\|\bar{x}\|}$$

$$= \frac{a^*(y^*)r}{2}.$$

This implies that

$$\theta^* \in \widehat{D}^*P_{r\mathbb{B}}(\bar{x})(y^*) \Rightarrow a^*(y^*) \leq 0. \quad (3.11)$$

Assume, by the way of contradiction, that $-J^*(o^*(y^*) \in \bar{x}_r^{\uparrow}$. Then there is $\delta > 0$ such that $\|\bar{x} - tJ^*(o^*(y^*)\| > r$, for all $t \in (0, \delta)$. In this case, we take a directional line segment in the limit (3.2) as $u = \bar{x} - tJ^*(o^*(y^*))$, for $t \downarrow 0$ with $t \in (0, \delta)$. By $\|\bar{x} - tJ^*(o^*(y^*))\| > r$, we have

$$P_{r\mathbb{B}}(\bar{x} - tJ^*(o^*(y^*))) = \frac{r}{\|\bar{x}-tJ^*(o^*(y^*))\|}(\bar{x} - tJ^*(o^*(y^*))), \text{ for } t \in (0, \delta).$$

So, by (3.10) and $\langle o^*(y^*), \bar{x}\rangle = 0$, it follows that

$$0 \geq \limsup_{u \to \bar{x}} \frac{-\langle y^*, P_{r\mathbb{B}}(u)-\bar{x}\rangle}{\|u-\bar{x}\|+\|P_{r\mathbb{B}}(u)-\bar{x}\|}$$

$$\geq \limsup_{t\downarrow 0, t<\delta} \frac{-\langle a^*(y^*)\bar{x}^* + o^*(y^*), P_{r\mathbb{B}}(\bar{x} - tJ^*(o^*(y^*))) - \bar{x}\rangle}{\|\bar{x} - tJ^*(o^*(y^*)) - \bar{x}\| + \|P_{r\mathbb{B}}(\bar{x} - tJ^*(o^*(y^*))) - \bar{x}\|}$$

$$= \limsup_{t\downarrow 0, t<\delta} \frac{-\langle a^*(y^*)\bar{x}^* + o^*(y^*),\ \frac{r}{\|\bar{x} - tJ^*(o^*(y^*))\|}(\bar{x} - tJ^*(o^*(y^*))) - \bar{x}\rangle}{t\|J^*(o^*(y^*))\| + \left\|\frac{r}{\|\bar{x} - tJ^*(o^*(y^*))\|}(\bar{x} - tJ^*(o^*(y^*))) - \bar{x}\right\|}$$

$$= \limsup_{t\downarrow 0, t<\delta} \frac{-a^*(y^*)\left(\frac{r}{\|\bar{x} - tJ^*(o^*(y^*))\|} - 1\right)\|\bar{x}\|^2 + \frac{rt\|J^*(o^*(y^*))\|^2}{\|\bar{x} - tJ^*(o^*(y^*))\|} + \frac{tra^*(y^*)}{\|\bar{x} - tJ^*(o^*(y^*))\|}\langle \bar{x}^*,\ J^*(o^*(y^*))\rangle}{t\|J^*(o^*(y^*))\| + \left\|\frac{r}{\|\bar{x} - tJ^*(o^*(y^*))\|}(\bar{x} - tJ^*(o^*(y^*))) - \bar{x}\right\|}$$

$$= \limsup_{t\downarrow 0, t<\delta} \frac{-a^*(y^*)r^2\left(\frac{r}{\|\bar{x} - tJ^*(o^*(y^*))\|} - 1\right) + \frac{rt\|J^*(o^*(y^*))\|^2}{\|\bar{x} - tJ^*(o^*(y^*))\|} + \frac{tra^*(y^*)}{\|\bar{x} - tJ^*(o^*(y^*))\|}\langle \bar{x}^*,\ J^*(o^*(y^*))\rangle}{t\|J^*(o^*(y^*))\| + \left\|\frac{r}{\|\bar{x} - tJ^*(o^*(y^*))\|}(\bar{x} - tJ^*(o^*(y^*))) - \bar{x}\right\|}$$

$$= \limsup_{t\downarrow 0, t<\delta} \frac{-a^*(y^*)r^2\left(\frac{r}{\|\bar{x} - tJ^*(o^*(y^*))\|} - \frac{r}{\|\bar{x}\|}\right) + \frac{tra^*(y^*)}{\|\bar{x} - tJ^*(o^*(y^*))\|}\langle \bar{x}^*,\ J^*(o^*(y^*))\rangle + \frac{rt\|J^*(o^*(y^*))\|^2}{\|\bar{x} - tJ^*(o^*(y^*))\|}}{t\|J^*(o^*(y^*))\| + \left\|\frac{r}{\|\bar{x} - tJ^*(o^*(y^*))\|}(\bar{x} - tJ^*(o^*(y^*))) - \bar{x}\right\|}$$

$$= \limsup_{t\downarrow 0, t<\delta} \frac{\frac{ta^*(y^*)r^3}{\|\bar{x} - tJ^*(o^*(y^*))\|\|\bar{x}\|}\left(\frac{\|\bar{x} + t(-J^*(o^*(y^*)))\| - \|\bar{x}\|}{t}\right) + \frac{tra^*(y^*)}{\|\bar{x} - tJ^*(o^*(y^*))\|}\langle \bar{x}^*,\ J^*(o^*(y^*))\rangle + \frac{tr\|J^*(o^*(y^*))\|^2}{\|\bar{x} - tJ^*(o^*(y^*))\|}}{t\|J^*(o^*(y^*))\| + \left\|\frac{r}{\|\bar{x} - tJ^*(o^*(y^*))\|}(\bar{x} - tJ^*(o^*(y^*))) - \bar{x}\right\|}$$

$$= \limsup_{t\downarrow 0, t<\delta} \frac{\frac{ta^*(y^*)r^3}{\|\bar{x} - tJ^*(o^*(y^*))\|\|\bar{x}\|}\left(\frac{\|\bar{x} + t(-J^*(o^*(y^*)))\| - \|\bar{x}\|}{t}\right) + \frac{tra^*(y^*)}{\|\bar{x} - tJ^*(o^*(y^*))\|}\langle \bar{x}^*,\ J^*(o^*(y^*))\rangle + \frac{tr\|J^*(o^*(y^*))\|^2}{\|\bar{x} - tJ^*(o^*(y^*))\|}}{t\|J^*(o^*(y^*))\| + \frac{tr}{\|\bar{x} - tJ^*(o^*(y^*))\|}\left\|\frac{\|\bar{x} + t(-J^*(o^*(y^*)))\| - \|\bar{x}\|}{t}\bar{x} + J^*(o^*(y^*))\right\|}$$

$$= \limsup_{t\downarrow 0, t<\delta} \frac{\frac{a^*(y^*)r^3}{\|\bar{x} - tJ^*(o^*(y^*))\|\|\bar{x}\|}\left(\frac{\|\bar{x} + t(-J^*(o^*(y^*)))\| - \|\bar{x}\|}{t}\right) + \frac{ra^*(y^*)}{\|\bar{x} - tJ^*(o^*(y^*))\|}\langle \bar{x}^*,\ J^*(o^*(y^*))\rangle + \frac{tr\|J^*(o^*(y^*))\|^2}{\|\bar{x} - tJ^*(o^*(y^*))\|}}{\|J^*(o^*(y^*))\| + \frac{r}{\|\bar{x} - tJ^*(o^*(y^*))\|}\left\|\frac{\|\bar{x} + t(-J^*(o^*(y^*)))\| - \|\bar{x}\|}{t}\bar{x} + J^*(o^*(y^*))\right\|}$$

$$= \frac{\frac{a^*(y^*)r^2}{\|\bar{x} - tJ^*(o^*(y^*))\|}\left(\frac{\langle J(\bar{x}),\ -J^*(o^*(y^*))\rangle}{\|\bar{x}\|}\right) + \frac{ra^*(y^*)}{\|\bar{x} - tJ^*(o^*(y^*))\|}\langle \bar{x}^*,\ J^*(o^*(y^*))\rangle + \frac{r\|J^*(o^*(y^*))\|^2}{\|\bar{x} - tJ^*(o^*(y^*))\|}}{\|J^*(o^*(y^*))\| + \left\|\left(\frac{\langle J(\bar{x}),\ -J^*(o^*(y^*))\rangle}{\|\bar{x}\|}\right)\bar{x} + J^*(o^*(y^*))\right\|}$$

$$= \frac{-\frac{a^*(y^*)r}{\|\bar{x}\|}\langle J(\bar{x}),\ J^*(o^*(y^*))\rangle + \frac{ra^*(y^*)}{\|\bar{x}\|}\langle \bar{x}^*,\ J^*(o^*(y^*))\rangle + \frac{r\|J^*(o^*(y^*))\|^2}{\|\bar{x}\|}}{\|J^*(o^*(y^*))\| + \left\|\left(\frac{\langle J(\bar{x}),\ -J^*(o^*(y^*))\rangle}{\|\bar{x}\|}\right)\bar{x} + J^*(o^*(y^*))\right\|}$$

$$= \frac{\|J^*(o^*(y^*))\|^2}{\|J^*(o^*(y^*))\| + \left\|\left(\frac{\langle J(\bar{x}),\ -J^*(o^*(y^*))\rangle}{\|\bar{x}\|}\right)\bar{x} + J^*(o^*(y^*))\right\|}$$

$$> 0.$$

This contradiction proves that

$$\theta^* \in \widehat{D}^*P_{r\mathbb{B}}(\bar{x})(y^*) \Longrightarrow -J^*(o^*(y^*)) \notin \bar{x}_r^\uparrow. \tag{3.12}$$

By (3.12), we must have $-J^*(o^*(y^*)) \in \bar{x}_r^\downarrow$. Then there is $\delta > 0$ such that $\|\bar{x} - tJ^*(o^*(y^*))\| \leq r$, for all $t \in (0, \delta)$. In this case, we take a directional line segment in the limit (3.2) as $u = \bar{x} - tJ^*(o^*(y^*))$, for $t \downarrow 0$ with $t \in (0, \delta)$. By $\|\bar{x} - tJ^*(o^*(y^*))\| \leq r$ in this case, we have

$$P_{r\mathbb{B}}(\bar{x} - tJ^*(o^*(y^*))) = \bar{x} - tJ^*(o^*(y^*)), \text{ for } t \in (0, \delta).$$

So, by (3.10) and $\langle o^*(y^*), \bar{x}\rangle = 0$, it follows that

$$0 \geq \limsup_{u \to \bar{x}} \frac{-\langle y^*, P_{r\mathbb{B}}(u)-\bar{x}\rangle}{\|u-\bar{x}\|+\|P_{r\mathbb{B}}(u)-\bar{x}\|}$$

$$\geq \limsup_{t \downarrow 0, t < \delta} \frac{-\langle a^*(y^*)\bar{x}^* + o^*(y^*), \bar{x} - tJ^*(o^*(y^*)) - \bar{x}\rangle}{\|\bar{x} - tJ^*(o^*(y^*)) - \bar{x}\| + \|\bar{x} - tJ^*(o^*(y^*)) - \bar{x}\|}$$

$$= \limsup_{t \downarrow 0, t < \delta} \frac{\langle a^*(y^*)\bar{x}^* + o^*(y^*), tJ^*(o^*(y^*))\rangle}{2t\|J^*(o^*(y^*))\|}$$

$$= \frac{a^*(y^*)\langle \bar{x}^*, J^*(o^*(y^*))\rangle + \|J^*(o^*(y^*))\|^2}{2\|J^*(o^*(y^*))\|}$$

$$= \frac{\frac{\langle y^*, \bar{x}\rangle}{\|\bar{x}\|^2}\langle \bar{x}^*, J^*(o^*(y^*))\rangle + \|J^*(o^*(y^*))\|^2}{2\|J^*(o^*(y^*))\|}$$

$$= \frac{\frac{\langle y^*, \bar{x}\rangle}{r^2}\langle \bar{x}^*, J^*(o^*(y^*))\rangle + \|J^*(o^*(y^*))\|^2}{2\|J^*(o^*(y^*))\|}.$$

This implies

$$\frac{\langle y^*, \bar{x}\rangle}{r^2}\langle \bar{x}^*, J^*(o^*(y^*))\rangle + \|J^*(o^*(y^*))\|^2 \leq 0. \tag{3.13}$$

Then, (3.9) is proved by (3.11). (3.12) and (3.13).

Next, we prove part 2 of (b) in (iii). That is, we prove that

$$o^*(y^*) = \theta^* \text{ and } \langle y^*, \bar{x}\rangle < 0 \implies \theta^* \in \widehat{D}^* P_{r\mathbb{B}}(\bar{x})(y^*).$$

Let $y^* = a^*(y^*)J(\bar{x})$ with $\langle y^*, \bar{x}\rangle < 0$, which implies $a^*(y^*) = \frac{\langle y^*, \bar{x}\rangle}{\|\bar{x}\|^2} < 0$. By $\|P_{r\mathbb{B}}(u)\| \leq r$, for any $u \in X$, we calculate

$$\limsup_{u \to \bar{x}} \frac{\langle(\theta^*, -y^*), (u, P_{r\mathbb{B}}(u)) - (\bar{x}, P_{r\mathbb{B}}(\bar{x}))\rangle}{\|u-\bar{x}\| + \|P_{r\mathbb{B}}(u) - P_{r\mathbb{B}}(\bar{x})\|}$$

$$= \limsup_{u \to \bar{x}} \frac{\langle(\theta^*, -a^*(y^*)J(\bar{x})), (u, P_{r\mathbb{B}}(u)) - (\bar{x}, P_{r\mathbb{B}}(\bar{x}))\rangle}{\|u-\bar{x}\| + \|P_{r\mathbb{B}}(u) - P_{r\mathbb{B}}(\bar{x})\|}$$

$$= \limsup_{u \to \bar{x}} \frac{\langle -a^*(y^*)J(\bar{x}), P_{r\mathbb{B}}(u) - \bar{x}\rangle}{\|u-\bar{x}\| + \|P_{r\mathbb{B}}(u) - \bar{x}\|}$$

$$= \limsup_{u \to \bar{x}} \frac{a^*(y^*)(\|\bar{x}\|^2 - \langle J(\bar{x}), P_{r\mathbb{B}}(u)\rangle)}{\|u-\bar{x}\| + \|P_{r\mathbb{B}}(u) - \bar{x}\|}$$

$$\leq \limsup_{u \to \bar{x}} \frac{a^*(y^*)(\|\bar{x}\|^2 - \|\bar{x}\|\|P_{r\mathbb{B}}(u)\|)}{\|u-\bar{x}\| + \|P_{r\mathbb{B}}(u) - \bar{x}\|}$$

$$\leq \limsup_{u \to \bar{x}} \frac{0}{\|u-\bar{x}\|+\| P_{r\mathbb{B}}(u)-\bar{x}\|}$$

$$\leq 0.$$

By definition, this implies $\theta^* \in \widehat{D}^* P_{r\mathbb{B}}(\bar{x})(y^*)$, for $y^* = \frac{\langle y^*,\bar{x}\rangle}{\|\bar{x}\|^2} J(\bar{x}) = a^*(y^*) J(\bar{x})$ with $a^*(y^*) < 0$, which proves part 2 of (b) in (iii). Then, part (b) in (iii) is proved.

Proof of (c). Take an arbitrary $z^* \in X^*$, by $\|\bar{x}\| = r$ and $J(\bar{x}) := \bar{x}^*$, we have

$$z^* \in \widehat{D}^* P_{r\mathbb{B}}(\bar{x})(J(\bar{x})) = \widehat{D}^* P_{r\mathbb{B}}(\bar{x})(\bar{x}^*) \iff \limsup_{u \to \bar{x}} \frac{\langle (z^*,-\bar{x}^*),\ (u,P_{r\mathbb{B}}(u)-(\bar{x},P_{r\mathbb{B}}(\bar{x})))\rangle}{\|u-\bar{x}\|+\| P_{r\mathbb{B}}(u)-P_{r\mathbb{B}}(\bar{x})\|} \leq 0$$

$$\iff \limsup_{u \to \bar{x}} \frac{\langle (z^*,-\bar{x}^*),\ (u,P_{r\mathbb{B}}(u)-(\bar{x},\bar{x}))\rangle}{\|u-\bar{x}\|+\| P_{r\mathbb{B}}(u)-\bar{x}\|} \leq 0$$

$$\iff \limsup_{u \to \bar{x}} \frac{\langle z^*,\ u-\bar{x}\rangle - \langle \bar{x}^*,\ P_{r\mathbb{B}}(u)-\bar{x}\rangle}{\|u-\bar{x}\|+\| P_{r\mathbb{B}}(u)-\bar{x}\|} \leq 0. \qquad (3.14)$$

We take a directional line segment in the limit in (3.14) $u = (1-t)\bar{x}$, for $t \downarrow 0$ with $t < 1$. Then, by Lemma 2.1, we have

$$0 \geq \limsup_{u \to \bar{x}} \frac{\langle z^*,\ u-\bar{x}\rangle - \langle \bar{x}^*,\ P_{r\mathbb{B}}(u)-\bar{x}\rangle}{\|u-\bar{x}\|+\| P_{r\mathbb{B}}(u)-\bar{x}\|}$$

$$\geq \limsup_{t \downarrow 0, t<1} \frac{\langle z^*,\ (1-t)\bar{x}-\bar{x}\rangle - \langle \bar{x}^*,\ (1-t)\bar{x}-\bar{x}\rangle}{\|(1-t)\bar{x}-\bar{x}\|+\|(1-t)\bar{x}-\bar{x}\|}$$

$$= \limsup_{t \downarrow 0, t<1} \frac{\langle a^*(z^*)\bar{x}^* + o^*(z^*), -t\bar{x}\rangle - \langle \bar{x}^*, -t\bar{x}\rangle}{2t\|\bar{x}\|}$$

$$= \frac{-a^*(z^*)\|\bar{x}\|^2 + \|\bar{x}\|^2}{2\|\bar{x}\|}$$

$$= \frac{(1-a^*(z^*))\|\bar{x}\|}{2}.$$

It implies that

$$z^* \in \widehat{D}^* P_{r\mathbb{B}}(\bar{x})(J(\bar{x})) \implies a^*(z^*) \geq 1. \qquad (3.15)$$

If we take a directional line segment in the limit in (3.14) as $u = (1+t)\bar{x}$, for $t \downarrow 0$. Then

$$0 \geq \limsup_{u \to \bar{x}} \frac{\langle z^*,\ u-\bar{x}\rangle - \langle \bar{x}^*,\ P_{r\mathbb{B}}(u)-\bar{x}\rangle}{\|u-\bar{x}\|+\| P_{r\mathbb{B}}(u)-\bar{x}\|}$$

$$\geq \limsup_{t \downarrow 0} \frac{\langle z^*,\ (1+t)\bar{x}-\bar{x}\rangle - \langle \bar{x}^*,\ \bar{x}-\bar{x}\rangle}{\|(1+t)\bar{x}-\bar{x}\|+\|\bar{x}-\bar{x}\|}$$

$$= \limsup_{t \downarrow 0} \frac{\langle a^*(z^*)\bar{x}^* + o^*(z^*),\ t\bar{x}\rangle - \langle \bar{x}^*, \theta\rangle}{t\|\bar{x}\|}$$

$$= \limsup_{t\downarrow 0} \frac{\langle a^*(z^*)\bar{x}^*, \ t\bar{x}\rangle + \langle o^*(z^*), \ t\bar{x}\rangle}{t\|\bar{x}\|}$$

$$= \frac{a^*(z^*)\|\bar{x}\|^2}{\|\bar{x}\|}$$

$$= a^*(z^*)r.$$

This implies that

$$z^* \in \widehat{D}^* P_{r\mathbb{B}}(\bar{x})(J(\bar{x})) \quad \Longrightarrow \quad a^*(z^*) \leq 0. \tag{3.16}$$

(3.16) contradicts to (3.15). This implies that $\widehat{D}^* P_{r\mathbb{B}}(\bar{x})(J(\bar{x})) = \emptyset$, which proves (c). □

In particular, if the underlying uniformly convex and uniformly smooth Banach space in Theorem 3.1 is a Hilbert space, then, it induces Theorem 3.2 in [14]. We show it below as a corollary of Theorem 3.1.

**Corollary 3.2** (Theorem 3.2 in [14]). *Let $H$ be a Hilbert space. For any $r > 0$, the Mordukhovich derivatives of $P_{r\mathbb{B}}$ at a point $\bar{x} \in H$ is given below.*

(I)      If $\bar{x} \in r\mathbb{B}^\circ$, then
$$\widehat{D}^* P_{r\mathbb{B}}(\bar{x})(y) = y, \text{ for every } y \in H.$$

(II)      If $\bar{x} \in H \setminus r\mathbb{B}$, then
$$\widehat{D}^* P_{r\mathbb{B}}(\bar{x})(y) = \frac{r}{\|\bar{x}\|}\left(y - \frac{\langle y,\bar{x}\rangle}{\|\bar{x}\|^2}\bar{x}\right), \text{ for every } y \in H.$$

(III)      If $\bar{x} \in r\mathbb{S}$, then

    (a)    $\widehat{D}^* P_{r\mathbb{B}}(\bar{x})(\theta) = \{\theta\};$

    (b)    *For any $y \in H \setminus \{\theta\}$, we have*
$$\theta \in \widehat{D}^* P_{r\mathbb{B}}(\bar{x})(y) \iff y = \frac{\langle y,\bar{x}\rangle}{\|\bar{x}\|^2}\bar{x} \text{ and } \langle y, \bar{x}\rangle < 0;$$

    (c)    $\widehat{D}^* P_{r\mathbb{B}}(\bar{x})(\bar{x}) = \emptyset.$

*Proof.* We only need to show that, if the underlying space is a Hilbert space, then part (b) in (iii) of Theorem 3.1 induces part (b) in (III) of Theorem 3.2 in [14]. Rest of the statements can be straight forwardly checked.

If $X$ is a Hilbert space in Theorem 3.1, then $\langle \bar{x}^*, \ J^*(o^*(y^*))\rangle = \langle \bar{x}, \ o(y)\rangle = 0$. From part 1 in (iii) of Theorem 3.1, this implies

$$\|o(y)\|^2 = \|J^*(o^*(y^*))\|^2 = 0.$$

That is, $o(y) = \theta$. Hence, part 1 in (iii) of Theorem 3.1 becomes

$$\theta \in \widehat{D}^* P_{r\mathbb{B}}(\bar{x})(y) \implies o(y) = \theta \text{ and } \langle y, \bar{x}\rangle < 0, \text{ for } y \in X \setminus \{\theta\}.$$

Then, by part 2 of in (iii) of Theorem 3.1, this implies

$$o(y) = \theta \text{ and } \langle y, \bar{x} \rangle < 0 \implies \theta \in \widehat{D}^* P_{r\mathbb{B}}(\bar{x})(y).$$

That is,

$$\theta \in \widehat{D}^* P_{r\mathbb{B}}(\bar{x})(y) \iff a(y) < 0 \text{ and } o(y) = \theta, \text{ for } y \in X \setminus \{\theta\}.$$

This coincides with part (III) of Theorem 3.2 in [14]. □

## 4. Mordukhovich derivatives of the metric projection onto closed and convex cylinders in real Banach space $l_p$

Let $(l_p, \|\cdot\|_p)$ be the real uniformly convex and uniformly smooth Banach space with dual space $(l_q, \|\cdot\|_q)$ satisfying $1 < p, q < \infty$ and $\frac{1}{p} + \frac{1}{q} = 1$. The origins of both $l_p$ and $l_q$ are exactly same $\theta = \theta^* = (0, 0, \dots)$. In this section, we recall some closed and convex cylinders in $l_p$ defined in [15], in which we studied the Fréchet differentiability of the metric projection onto closed and convex cylinders in $l_p$. In this paper, we use the results from [15] to investigate Mordukhovich derivatives of the metric projection onto closed and convex cylinders in $l_p$. At first, we recall some concepts and notations used in [15].

Let $J: l_p \to l_q$ be the normalized duality mapping in $l_p$. It has the following representations. For any point $x = (x_1, x_2, \dots) \in l_p$ with $x \neq \theta$, we have

$$J(x) = \left( \frac{|x_1|^{p-1} \text{sign}(x_1)}{\|x\|_p^{p-2}}, \frac{|x_2|^{p-1} \text{sign}(x_2)}{\|x\|_p^{p-2}}, \dots \right)$$

$$= \left( \frac{|x_1|^{p-2} x_1}{\|x\|_p^{p-2}}, \frac{|x_2|^{p-2} x_2}{\|x\|_p^{p-2}}, \dots \right). \tag{4.1}$$

Similarly, to (4.1), the representations of the normalized duality mapping $J^*: l_q \to l_p$ is given, for any $y^* = (y_1, y_2, \dots) \in l_q$ with $y \neq \theta$, by

$$J^*(y^*) = \left( \frac{|y_1|^{q-1} \text{sign}(y_1)}{\|y\|_q^{q-2}}, \frac{|y_2|^{q-1} \text{sign}(y_2)}{\|y\|_q^{q-2}}, \dots \right)$$

$$= \left( \frac{|y_1|^{q-2} y_1}{\|y\|_q^{q-2}}, \frac{|y_2|^{q-2} y_2}{\|y\|_q^{q-2}}, \dots \right). \tag{4.2}$$

Let $\mathbb{N}$ denote the set of all positive integers. Let $M$ be a nonempty subset of $\mathbb{N}$ and let $\bar{M} = \mathbb{N} \setminus M$, which is the complementary set of $M$. Following [15], we define

$$l_p^M = \{x = (x_1, x_2, \dots) \in l_p : x_i = 0, \text{ for all } i \in \bar{M}\},$$

$$l_q^M = \{y = (y_1, y_2, \dots) \in l_q : y_i = 0, \text{ for all } i \in \bar{M}\}.$$

$l_p^M$ and $l_q^M$ are closed subspaces of $l_p$ and $l_q$, respectively. They are the dual spaces of each other.

We define a mapping $(\cdot)_M: l_p \to l_p^M$, for $x = (x_1, x_2, \ldots) \in l_p$, by

$$(x_M)_i = \begin{cases} x_i, & \text{for } i \in M, \\ 0, & \text{for } i \notin M, \end{cases} \quad \text{for } i \in \mathbb{N}.$$

Similarly, we define a mapping $(\cdot)_{\bar{M}}: l_p \to l_p^{\bar{M}}$, for $x = (x_1, x_2, \ldots) \in l_p$ by

$$(x_{\bar{M}})_i = \begin{cases} x_i, & \text{for } i \in \bar{M}, \\ 0, & \text{for } i \notin \bar{M}, \end{cases} \quad \text{for } i \in \mathbb{N}.$$

Then, $l_p$ and $l_q$ have the following decompositions

$$x = x_M + x_{\bar{M}}, \text{ for any } x \in l_p.$$

and
$$x^* = x_M^* + x_{\bar{M}}^*, \text{ for any } x^* \in l_q. \tag{4.3}$$

**Lemma 4.1 in [15]**. *Let $M$ be a nonempty subset of $\mathbb{N}$. Then $J$ is the normalized duality mapping from $l_p^M$ to $l_q^M$. That is,*

$$J(x) \in l_q^M, \text{ for any } x \in l_p^M.$$

Let $\mathbb{B}_M$ denote the unit closed ball in $l_p^M$. It follows that, for any $r > 0$, $r\mathbb{B}_M$ is the closed ball with radius $r$ and center origin in $l_p^M$. Let $\mathbb{S}_M$ be the unit sphere in $l_p^M$. Then, $r\mathbb{S}_M$ is the sphere in $l_p^M$ with radius $r$ and center $\theta$. We define

$$\mathbb{C}_M = \{x \in l_p: x_M \in \mathbb{B}_M\}.$$

$\mathbb{C}_M$ is called the cylinder in $l_p$ with base $\mathbb{B}_M$. It is a closed and convex subset in $l_p$. For any $r > 0$, $r\mathbb{C}_M$ is the cylinder in $l_p$ with base $r\mathbb{B}_M$, which is a closed and convex subset of $l_p$. More precisely speaking, we have

$$r\mathbb{C}_M = \{x \in l_p: x_M \in r\mathbb{B}_M\}.$$

This implies that, for any $x \in l_p$, we have

$$x \in r\mathbb{C}_M \quad \Leftrightarrow \quad \|x_M\|_p \leq r. \tag{4.4}$$

The boundary of $r\mathbb{C}_M$ is denoted by $\partial(r\mathbb{C}_M)$ satisfying

$$\partial(r\mathbb{C}_M) = \{x \in l_p: \|x_M\|_p = r\}.$$

**Lemma 4.2 in [15]**. *For any $r > 0$, the metric projection $P_{r\mathbb{C}_M}: l_p \to r\mathbb{C}_M$ satisfies the following formula.*

$$P_{r\mathbb{C}_M}(x) = \begin{cases} x, & \text{for any } x \in r\mathbb{C}_M, \\ \dfrac{r}{\|x_M\|_p} x_M + x_{\bar{M}}, & \text{for any } x \in l_p \setminus r\mathbb{C}_M. \end{cases} \tag{4.5}$$

For any $x \in l_p$ with $\|x_M\|_p = r$, two subsets $x_r^{\Uparrow}$ and $x_r^{\Downarrow}$ of $l_p$ are defined by

(a) $x_r^{\Uparrow} = \{v \in l_p:$ there is $\delta > 0$ such that $\|(x+tv)_M\|_p > r$, for all $t \in (0, \delta)\}$;
(b) $x_r^{\Downarrow} = \{v \in l_p:$ there is $\delta > 0$ such that $\|(x+tv)_M\|_p \leq r$, for all $t \in (0, \delta)\}$.

**Theorem 4.3 (partial) in [15].** *For any $r > 0$, the metric projection $P_{r\mathbb{C}_M}: l_p \to r\mathbb{C}_M$ has the following differentiable properties.*

(i) $P_{r\mathbb{C}_M}$ *is strictly Fréchet differentiable on $(r\mathbb{C}_M)^\circ$ satisfying*

$$\nabla P_{r\mathbb{C}_M}(\bar{x}) = I_{l_p}, \text{ for every } \bar{x} \in (r\mathbb{C}_M)^\circ.$$

*That is,*

$$\bar{x} \in (r\mathbb{C}_M)^\circ \implies \nabla P_{r\mathbb{C}_M}(\bar{x})(u) = u, \text{ for every } u \in l_p.$$

(ii) $P_{r\mathbb{C}_M}$ *is Fréchet differentiable at every point $\bar{x} \in l_p \setminus r\mathbb{C}_M$ such that*

$$\nabla P_{r\mathbb{C}_M}(\bar{x})(u) = \frac{r}{\|\bar{x}_M\|_p}\left(u_M - \frac{\langle J(\bar{x}_M), u_M\rangle}{\|\bar{x}_M\|_p^2}\bar{x}_M\right) + u_{\bar{M}}, \text{ for every } u \in l_p. \qquad (4.6)$$

*In particular,*

$$\nabla P_{r\mathbb{C}_M}(\bar{x})(\bar{x}) = x_{\bar{M}}, \text{ for every } \bar{x} \in l_p \setminus r\mathbb{C}_M.$$

(iii) $P_{r\mathbb{C}_M}$ *is not Fréchet differentiable at any point $\bar{x} \in \partial(r\mathbb{C}_M)$, that is,*

$$\nabla P_{r\mathbb{C}_M}(\bar{x}) \text{ does not exist, for any } \bar{x} \in \partial(r\mathbb{C}_M).$$

By Theorem 4.3 in [15] restated as above, we investigate the Mordukhovich derivatives of the metric projection $P_{r\mathbb{C}_M}$ onto closed and convex cylinders $r\mathbb{C}_M$ in $l_p$.

**Theorem 4.1.** *For any $r > 0$, Mordukhovich derivatives of the metric projection $P_{r\mathbb{C}_M}: l_p \to r\mathbb{C}_M$ have the following representations.*

(i) *For every $\bar{x} \in (r\mathbb{C}_M)^\circ$, we have*

$$\widehat{D}^* P_{r\mathbb{C}_M}(\bar{x})(y^*) = \{y^*\}, \text{ for every } y^* \in l_q.$$

(ii) *For every point $\bar{x} \in l_p \setminus r\mathbb{C}_M$, we have*

$$\widehat{D}^* P_{r\mathbb{C}_M}(\bar{x})(y^*) = \left\{\frac{r}{\|\bar{x}_M\|_p}\left(y_M^* - \frac{\langle y_M^*, \bar{x}_M\rangle}{\|\bar{x}_M\|_p^2}J(\bar{x}_M)\right) + y_{\bar{M}}^*\right\}, \text{ for all } y^* \in l_q.$$

*In particular,*

$$\widehat{D}^* P_{r\mathbb{C}_M}(\bar{x})(J(\bar{x})) = J(\bar{x})_{\bar{M}}, \text{ for every } \bar{x} \in l_p \setminus r\mathbb{C}_M.$$

(iii) *For any $\bar{x} \in \partial(r\mathbb{C}_M)$, we have*

(a) $\widehat{D}^* P_{r\mathbb{C}_M}(\bar{x})(\theta^*) = \{\theta^*\}$;

(b) *For any $y^* \in l_q$,*

$$\theta^* \in \widehat{D}^* P_{r\mathbb{C}_M}(\bar{x})(y^*)$$

$$\Leftrightarrow \quad y_{\bar{M}}^* = \theta^*, \ -(J^*(y^*))_M \notin \bar{x}_r^{\Downarrow} \ \text{and} \ \langle y_M^*, \bar{x}_M \rangle = -r\|y_M^*\|_q.$$

*In particular,*

$$\theta^* \in \widehat{D}^* P_{r\mathbb{C}_M}(\bar{x})(-J(\bar{x})_M);$$

(c) $\widehat{D}^* P_{r\mathbb{C}_M}(\bar{x})(\bar{x}) = \emptyset$.

*Proof.* Proof of (i). By part (i) in Theorem 4.3 in [15], we have

$$\nabla P_{r\mathbb{C}_M}(\bar{x}) = I_{l_p}, \text{ for every } \bar{x} \in (r\mathbb{C}_M)^\circ.$$

Then, by Theorem 1.38 in [18], for every $\bar{x} \in (r\mathbb{C}_M)^\circ$, we have

$$\widehat{D}^* P_{r\mathbb{C}_M}(\bar{x})(y^*) = \{(I_{l_p})^*(y^*)\} = \{I_{l_q}(y^*)\} = \{y^*\}, \text{ for all } y^* \in l_q.$$

Proof of (ii). By part (ii) in Theorem 4.3 in [15], for every $\bar{x} \in l_p \setminus r\mathbb{C}_M$, we have

$$\nabla P_{r\mathbb{C}_M}(\bar{x})(u) = \frac{r}{\|\bar{x}_M\|_p}\left(u_M - \frac{\langle J(\bar{x}_M), u_M \rangle}{\|\bar{x}_M\|_p^2}\bar{x}_M\right) + u_{\bar{M}}, \text{ for every } u \in l_p. \qquad (4.6)$$

Then, by Theorem 1.38 in [18], for every $\bar{x} \in l_p \setminus r\mathbb{C}_M$, we have

$$\widehat{D}^* P_{r\mathbb{C}_M}(\bar{x})(y^*) = \{(\nabla P_{r\mathbb{C}_M}(\bar{x}))^*(y^*)\}, \text{ for all } y^* \in l_q.$$

Since $\left(\nabla P_{r\mathbb{C}_M}(\bar{x})\right)^*(y^*) \in l_q$, then, for any $u \in l_p$, by (4.6) and by the decompositions (4.3) and $\langle J(\bar{x}_M), u_{\bar{M}} \rangle = 0$, we have

$$\langle \left(\nabla P_{r\mathbb{C}_M}(\bar{x})\right)^*(y^*), u \rangle$$

$$= \langle y^*, \nabla P_{r\mathbb{C}_M}(\bar{x})(u) \rangle$$

$$= \langle y^*, \frac{r}{\|\bar{x}_M\|_p}\left(u_M - \frac{\langle J(\bar{x}_M), u_M \rangle}{\|\bar{x}_M\|_p^2}\bar{x}_M\right) + u_{\bar{M}} \rangle$$

$$= \frac{r}{\|\bar{x}_M\|_p}\left(\langle y^*, u_M \rangle - \frac{\langle J(\bar{x}_M), u_M \rangle}{\|\bar{x}_M\|_p^2}\langle y^*, \bar{x}_M \rangle\right) + \langle y^*, u_{\bar{M}} \rangle$$

$$= \frac{r}{\|\bar{x}_M\|_p}\left(\langle y^*, u_M \rangle - \frac{\langle J(\bar{x}_M), u_M + u_{\bar{M}} \rangle}{\|\bar{x}_M\|_p^2}\langle y^*, \bar{x}_M \rangle\right) + \langle y^*, u_{\bar{M}} \rangle$$

$$= \frac{r}{\|\bar{x}_M\|_p}\left(\langle y^*, u_M\rangle - \frac{\langle y^*, \bar{x}_M\rangle}{\|\bar{x}_M\|_p^2}\langle J(\bar{x}_M), u\rangle\right) + \langle y^*, u_{\bar{M}}\rangle$$

$$= \frac{r}{\|\bar{x}_M\|_p}\left(\langle y_M^*, u_M\rangle - \frac{\langle y_M^*, \bar{x}_M\rangle}{\|\bar{x}_M\|_p^2}\langle J(\bar{x}_M), u\rangle\right) + \langle y_{\bar{M}}^*, u_{\bar{M}}\rangle$$

$$= \frac{r}{\|\bar{x}_M\|_p}\left(\langle y_M^*, u\rangle - \frac{\langle y_M^*, \bar{x}_M\rangle}{\|\bar{x}_M\|_p^2}\langle J(\bar{x}_M), u\rangle\right) + \langle y_{\bar{M}}^*, u\rangle$$

$$= \frac{r}{\|\bar{x}_M\|_p}\langle y_M^* - \frac{\langle y_M^*, \bar{x}_M\rangle}{\|\bar{x}_M\|_p^2}J(\bar{x}_M), u\rangle + \langle y_{\bar{M}}^*, u\rangle$$

$$= \langle \frac{r}{\|\bar{x}_M\|_p}\left(y_M^* - \frac{\langle y_M^*, \bar{x}_M\rangle}{\|\bar{x}_M\|_p^2}J(\bar{x}_M)\right) + y_{\bar{M}}^*, u\rangle, \text{ for any } u \in l_p.$$

This implies

$$\widehat{D}^* P_{r\mathbb{C}_M}(\bar{x})(y^*) = \left\{\left(\nabla P_{r\mathbb{C}_M}(\bar{x})\right)^*(y^*)\right\}$$

$$= \left\{\frac{r}{\|\bar{x}_M\|_p}\left(y_M^* - \frac{\langle y_M^*, \bar{x}_M\rangle}{\|\bar{x}_M\|_p^2}J(\bar{x}_M)\right) + y_{\bar{M}}^*\right\}, \text{ for all } y^* \in l_q.$$

Part (ii) is proved.

Proof of (a) in part (iii). For any $\bar{x} \in \partial(r\mathbb{C}_M)$, by $\|\bar{x}_M\|_p = r$, we have $P_{r\mathbb{C}_M}(\bar{x}) = \bar{x}$. It is clear that

$$\theta^* \in \widehat{D}^* P_{r\mathbb{C}_M}(\bar{x})(\theta^*). \tag{4.7}$$

Next, we prove $z^* \notin \widehat{D}^* P_{r\mathbb{C}_M}(\bar{x})(\theta^*)$, for $z^* \in l_q\setminus\{\theta^*\}$. By (2.20), we have

$$z^* \in \widehat{D}^* P_{r\mathbb{C}_M}(\bar{x})(\theta^*) \iff \limsup_{u\to\bar{x}} \frac{\langle(z^*,-\theta^*),\ (u, P_{r\mathbb{C}_M}(u) - (\bar{x}, P_{r\mathbb{C}_M}(\bar{x})))\rangle}{\|u-\bar{x}\|_p + \|P_{r\mathbb{C}_M}(u) - P_{r\mathbb{C}_M}(\bar{x})\|_p} \leq 0$$

$$\iff \limsup_{u\to\bar{x}} \frac{\langle(z^*,-\theta^*),\ (u, P_{r\mathbb{C}_M}(u) - (\bar{x}, P_{r\mathbb{C}_M}(\bar{x})))\rangle}{\|u-\bar{x}\|_p + \|P_{r\mathbb{C}_M}(u) - \bar{x}\|_p} \leq 0$$

$$\iff \limsup_{u\to\bar{x}} \frac{\langle z^*,\ u-\bar{x}\rangle}{\|u-\bar{x}\|_p + \|P_{r\mathbb{C}_M}(u) - \bar{x}\|_p} \leq 0. \tag{4.8}$$

Let $z := J^*(z^*)$. Since $z^* \in l_q\setminus\{\theta^*\}$, then

$$\langle z^*,\ J^*(z^*)\rangle = \langle z^*,\ z\rangle = \|z\|_p^2 = \|z^*\|_q^2 > 0.$$

Case 1. Suppose $z \in \bar{x}_r^\Uparrow$. Then, there is $\delta > 0$ such that $\|(x + tv)_M\|_p > r$, for all $t \in (0, \delta)$. In this case, we take a directional line segment in the limit (4.8) as $u = \bar{x} + tz$, for $t \downarrow 0$ with $t < \delta$. By (4.5) and $\|\bar{x}_M\|_p = r$, we have

$$P_{r\mathbb{C}_M}(\bar{x}+tz) = \frac{r}{\|(\bar{x}+tz)_M\|_p}(\bar{x}+tz)_M + (\bar{x}+tz)_{\bar{M}}$$

$$= \frac{r}{\|(\bar{x}+tz)_M\|_p}\bar{x}_M + \bar{x}_{\bar{M}} + \frac{tr}{\|(\bar{x}+tz)_M\|_p}z_M + tz_{\bar{M}}.$$

By (4.8), this implies

$$z^* \in \widehat{D}^*P_{r\mathbb{C}_M}(\bar{x})(\theta^*) \iff 0 \geq \limsup_{u \to \bar{x}} \frac{\langle z^*, u-\bar{x}\rangle}{\|u-\bar{x}\|_p + \|P_{r\mathbb{C}_M}(u)-\bar{x}\|_p}$$

$$\geq \limsup_{t\downarrow 0, t<\delta} \frac{\langle z^*, \bar{x}+tz-\bar{x}\rangle}{\|\bar{x}+tz-\bar{x}\|_p + \|P_{r\mathbb{C}_M}(\bar{x}+tz)-\bar{x}\|_p}$$

$$= \limsup_{t\downarrow 0, t<\delta} \frac{t\langle z^*, z\rangle}{t\|z\|_p + \left\|\frac{r}{\|\bar{x}+tz\|}(\bar{x}+tz)-\bar{x}\right\|_p}$$

$$= \limsup_{t\downarrow 0, t<\delta} \frac{t\langle z^*, z\rangle}{t\|z\|_p + \left\|\frac{r}{\|(\bar{x}+tz)_M\|_p}\bar{x}_M + \bar{x}_{\bar{M}} + \frac{tr}{\|(\bar{x}+tz)_M\|_p}z_M + tz_{\bar{M}} - \bar{x}\right\|_p}$$

$$= \limsup_{t\downarrow 0, t<\delta} \frac{t\langle z^*, z\rangle}{t\|z\|_p + \left\|\frac{r}{\|(\bar{x}+tz)_M\|_p}\bar{x}_M - \bar{x}_M + \frac{tr}{\|(\bar{x}+tz)_M\|_p}z_M + tz_{\bar{M}}\right\|_p}. \quad (4.9)$$

By $\|\bar{x}+tz\| > r$, for all $0 < t < \delta$, we calculate $\left\|\frac{r}{\|(\bar{x}+tz)_M\|_p}\bar{x}_M - \bar{x}_M + \frac{tr}{\|(\bar{x}+tz)_M\|_p}z_M + tz_{\bar{M}}\right\|_p$ in the denominator in limit (4.9).

$$\left\|\frac{r}{\|(\bar{x}+tz)_M\|_p}\bar{x}_M - \bar{x}_M + \frac{tr}{\|(\bar{x}+tz)_M\|_p}z_M + tz_{\bar{M}}\right\|_p$$

$$\leq \|\bar{x}_M\|_p \left|\frac{r}{\|(\bar{x}+tz)_M\|_p} - 1\right| + t\left\|\frac{r}{\|(\bar{x}+tz)_M\|_p}z_M + z_{\bar{M}}\right\|_p$$

$$= r\left|\frac{\|(\bar{x}+tz)_M\|_p - r}{\|(\bar{x}+tz)_M\|_p}\right| + t\left\|\frac{r}{\|(\bar{x}+tz)_M\|_p}z_M + z_{\bar{M}}\right\|_p$$

$$= r\left|\frac{\|\bar{x}_M+tz_M\|_p - \|\bar{x}_M\|_p}{\|(\bar{x}+tz)_M\|_p}\right| + t\left\|\frac{r}{\|(\bar{x}+\delta z)_M\|_p}z_M + z_{\bar{M}}\right\|_p$$

$$= \frac{tr}{\|(\bar{x}+tz)_M\|_p}\left|\frac{\|\bar{x}_M+tz_M\|_p - \|\bar{x}_M\|_p}{t}\right| + t\left\|\frac{r}{\|(\bar{x}+tz)_M\|_p}z_M + z_{\bar{M}}\right\|_p. \quad (4.10)$$

Substituting (4.10) into the limit (4.9), we have

$$z^* \in \widehat{D}^*P_{r\mathbb{C}_M}(\bar{x})(\theta^*) \iff 0 \geq \limsup_{t\downarrow 0, t<\delta} \frac{t\langle z^*, z\rangle}{t\|z\|_p + \left\|\frac{r}{\|(\bar{x}+tz)_M\|_p}\bar{x}_M - \bar{x}_M + \frac{tr}{\|(\bar{x}+tz)_M\|_p}z_M + tz_{\bar{M}}\right\|_p}$$

$$\geq \limsup_{t\downarrow 0, t<\delta} \frac{t\|z\|_p^2}{t\|z\|_p + \frac{tr}{\|(\bar{x}+tz)_M\|_p}\left|\frac{\|\bar{x}_M+tz_M\|_p - \|\bar{x}_M\|_p}{t}\right| + t\left\|\frac{r}{\|(\bar{x}+tz)_M\|_p} z_M + z_{\overline{M}}\right\|_p}$$

$$\geq \limsup_{t\downarrow 0, t<\delta} \frac{\|z\|_p^2}{\|z\|_p + \frac{r}{\|(\bar{x}+tz)_M\|_p}\left|\frac{\|\bar{x}_M+tz_M\|_p - \|\bar{x}_M\|_p}{t}\right| + \left\|\frac{r}{\|(\bar{x}+tz)_M\|_p} z_M + z_{\overline{M}}\right\|_p}$$

$$= \frac{\|z\|_p^2}{\|z\|_p + \left|\frac{\langle J(\bar{x}_M), z_M\rangle}{\|\bar{x}_M\|_p}\right| + \|z\|_p}$$

$$= \frac{\|z\|_p^2}{2\|z\|_p + \left|\frac{\langle J(\bar{x}_M), z_M\rangle}{r}\right|}$$

$$\geq \frac{\|z\|_p^2}{2\|z\|_p + \frac{1}{r}\|J(\bar{x}_M)\|_q \|z\|_p}$$

$$= \frac{\|z\|_p^2}{2\|z\|_p + \frac{1}{r}\|\bar{x}_M\|_p \|z\|_p}$$

$$= \frac{\|z\|_p}{3}$$

$$> 0. \tag{4.11}$$

(4.11) is a contradiction.

Case 2. Suppose $z \in x_r^{\Downarrow}$. Then, there is $\delta > 0$ such that $\|(x+tv)_M\|_p \leq r$, for all $t \in (0, \delta)$. In this case, we take a directional line segment in the limit (4.8) as $u = \bar{x} + tz$, for $t \downarrow 0$ with $t < \delta$. By (4.5) and by $\|\bar{x}_M\|_p = r$, we have

$$P_{r\mathbb{C}_M}(\bar{x} + tz) = \bar{x} + tz, \text{ for all } t \in (0, \delta).$$

By (4.8), this implies

$$z^* \in \widehat{D}^* P_{r\mathbb{C}_M}(\bar{x})(\theta^*) \iff 0 \geq \limsup_{u \to \bar{x}} \frac{\langle z^*, u-\bar{x}\rangle}{\|u-\bar{x}\|_p + \|P_{r\mathbb{C}_M}(u)-\bar{x}\|_p}$$

$$\geq \limsup_{t\downarrow 0, t<\delta} \frac{\langle z^*, \bar{x}+tz-\bar{x}\rangle}{\|\bar{x}+tz-\bar{x}\|_p + \|P_{r\mathbb{C}_M}(\bar{x}+tz)-\bar{x}\|_p}$$

$$= \limsup_{t\downarrow 0, t<\delta} \frac{t\langle z^*, z\rangle}{t\|z\|_p + \|\bar{x}+tz-\bar{x}\|_p}$$

$$= \limsup_{t\downarrow 0, t<\delta} \frac{t\|z\|_p^2}{2t\|z\|_p}$$

$$= \frac{\|z\|_p}{2}$$

$$> 0. \tag{4.12}$$

(4.12) is a contradiction. By (4.11) and (4.12), we have

$$z^* \in l_q \setminus \{\theta^*\} \implies z^* \notin \widehat{D}^* P_{r\mathbb{B}}(\bar{x})(\theta^*). \tag{4.13}$$

By (4.7) and (4.13), we proved (a) in (iii):

$$\widehat{D}^* P_{r\mathbb{C}_M}(\bar{x})(\theta^*) = \{\theta^*\}.$$

Proof of (b) in (iii). For any $y^* \in l_q$, from part (a) of (iii) in this theorem, we have

$$y^* = \theta^* \implies \theta^* \in \widehat{D}^* P_{r\mathbb{B}}(\bar{x})(y^*). \tag{4.14}$$

Then, we only need to consider $y^* \in l_q \setminus \{\theta^*\}$. We firstly prove

$$\theta^* \in \widehat{D}^* P_{r\mathbb{B}}(\bar{x})(y^*)$$
$$\implies y_{\bar{M}}^* = \theta^*, -(J^*(y^*))_M \notin \bar{x}_r^{\Downarrow} \text{ and } \langle y_M^*, \bar{x}_M \rangle = -\|y_M^*\|_q \|\bar{x}_M\|_p. \tag{4.15}$$

By (2.20), we have

$$\theta^* \in \widehat{D}^* P_{r\mathbb{C}_M}(\bar{x})(y^*) \iff \limsup_{u \to \bar{x}} \frac{\langle (\theta^*, -y^*), (u, P_{r\mathbb{C}_M}(u) - (\bar{x}, P_{r\mathbb{C}_M}(\bar{x}))) \rangle}{\|u - \bar{x}\|_p + \|P_{r\mathbb{C}_M}(u) - P_{r\mathbb{C}_M}(\bar{x})\|_p} \leq 0$$

$$\iff \limsup_{u \to \bar{x}} \frac{-\langle y^*, P_{r\mathbb{C}_M}(u) - P_{r\mathbb{C}_M}(\bar{x}) \rangle}{\|u - \bar{x}\|_p + \|P_{r\mathbb{C}_M}(u) - \bar{x}\|_p} \leq 0$$

$$\iff \limsup_{u \to \bar{x}} \frac{-\langle y^*, P_{r\mathbb{C}_M}(u) - \bar{x} \rangle}{\|u - \bar{x}\|_p + \|P_{r\mathbb{C}_M}(u) - \bar{x}\|_p} \leq 0. \tag{4.16}$$

At first, we prove

$$\theta^* \in \widehat{D}^* P_{r\mathbb{B}}(\bar{x})(y^*) \implies y_{\bar{M}}^* = \theta^*. \tag{4.15.1}$$

Assume, by the way of contradiction, that $\theta^* \in \widehat{D}^* P_{r\mathbb{B}}(\bar{x})(y^*)$ and $y_{\bar{M}}^* \neq \theta^*$. Notice that

$$\left\|(\bar{x} - tJ^*(y_{\bar{M}}^*))_M\right\|_p = \|\bar{x}_M\|_p = r, \text{ for all } t > 0. \tag{4.17}$$

Then, we take a directional line segment in the limit (4.16) as $u = \bar{x} - ty_{\bar{M}}^*$, for $t \downarrow 0$ with $t < 1$. By (4.17), $\bar{x} - ty_{\bar{M}}^* \in r\mathbb{C}_M$, for all $0 < t < 1$. we have

$$0 \geq \limsup_{u \to \bar{x}} \frac{-\langle y^*, P_{r\mathbb{C}_M}(u) - \bar{x} \rangle}{\|u - \bar{x}\|_p + \|P_{r\mathbb{C}_M}(u) - \bar{x}\|_p}$$

$$\geq \limsup_{t \downarrow 0, t < 1} \frac{-\langle y^*, \bar{x} - tJ^*(y_{\bar{M}}^*) - \bar{x} \rangle}{\|\bar{x} + tJ^*(y_{\bar{M}}^*) - \bar{x}\|_p + \|\bar{x} - tJ^*(y_{\bar{M}}^*) - \bar{x}\|_p}$$

$$= \limsup_{t\downarrow 0, t<1} \frac{t\langle y_{\overline{M}}^*, J^*(y_{\overline{M}}^*)\rangle}{2t\|J^*(y_{\overline{M}}^*)\|_p}$$

$$= \frac{\langle y_{\overline{M}}^*, J^*(y_{\overline{M}}^*)\rangle}{2\|J^*(y_{\overline{M}}^*)\|_p}$$

$$= \frac{\|J^*(y_{\overline{M}}^*)\|_p}{2}$$

$$> 0.$$

This contradiction proves (4.15.1), which is

$$\theta^* \in \widehat{D}^* P_{r\mathbb{B}}(\bar{x})(y^*) \Longrightarrow y_{\overline{M}}^* = \theta^*.$$

Since $y^* \in l_q\setminus\{\theta^*\}$ and $y_{\overline{M}}^* = \theta^*$, then, we must have $y_M^* \neq \theta^*$. This implies $J^*(y_M^*) \neq \theta$, which implies $(J^*(y^*))_M \neq \theta$. Next, we prove

$$\theta^* \in \widehat{D}^* P_{r\mathbb{C}_M}(\bar{x})(y^*) \Longrightarrow -(J^*(y^*))_M \notin \bar{x}_r^{\Downarrow}. \tag{4.15.2}$$

Assume, by the way of contradiction, that $-(J^*(y^*))_M \in \bar{x}_r^{\Downarrow}$. Then, there is $\delta > 0$ such that $\|\bar{x} - t(J^*(y^*))_M\|_p \leq r$, for all $t \in (0, \delta)$. In this case, we take a directional line segment in the limit (4.16) as $u = \bar{x} - t(J^*(y^*))_M$, for $t \downarrow 0$ with $t < \delta$. By (4.5) we have

$$P_{r\mathbb{C}_M}(\bar{x} - t(J^*(y^*))_M) = \bar{x} - t(J^*(y^*))_M, \text{ for all } 0 < t < \delta.$$

By $\|\bar{x}_M\|_p = r$, we have

$$0 \geq \limsup_{u\to\bar{x}} \frac{-\langle y^*, P_{r\mathbb{C}_M}(u) - \bar{x}\rangle}{\|u-\bar{x}\|_p + \|P_{r\mathbb{C}_M}(u)-\bar{x}\|_p}$$

$$\geq \limsup_{t\downarrow 0, t<\delta} \frac{-\langle y^*, \bar{x}-t(J^*(y^*))_M - \bar{x}\rangle}{\|\bar{x}-t(J^*(y^*))_M-\bar{x}\|_p + \|\bar{x}-t(J^*(y^*))_M-\bar{x}\|_p}$$

$$= \limsup_{t\downarrow 0, t<\delta} \frac{t\langle y^*, (J^*(y^*))_M\rangle}{2t\|(J^*(y^*))_M\|_p}$$

$$= \frac{\|(J^*(y^*))_M\|_p}{2}$$

$$> 0.$$

This contradiction proves (4.15.2). To the end of the proof of (4.15), we finally prove

$$\theta^* \in \widehat{D}^* P_{r\mathbb{C}_M}(\bar{x})(y^*) \Longrightarrow \langle y_M^*, \bar{x}_M\rangle = -\|y_M^*\|_q \|\bar{x}_M\|_p. \tag{4.15.3}$$

We take a directional line segment in the limit (4.16) as $u = \bar{x} - t\bar{x}_M$, for $t \downarrow 0$ with $t < 1$. By

$\bar{x} - t\bar{x}_M \in r\mathbb{C}_M$, for all $0 < t < 1$, we have

$$0 \geq \limsup_{u \to \bar{x}} \frac{-\langle y^*, P_{r\mathbb{C}_M}(u) - \bar{x}\rangle}{\|u - \bar{x}\|_p + \|P_{r\mathbb{C}_M}(u) - \bar{x}\|_p}$$

$$\geq \limsup_{t \downarrow 0, t < 1} \frac{-\langle y^*, \bar{x} - t\bar{x}_M - \bar{x}\rangle}{\|\bar{x} + t\bar{x}_M - \bar{x}\|_p + \|\bar{x} - t\bar{x}_M - \bar{x}\|_p}$$

$$= \limsup_{t \downarrow 0, t < 1} \frac{t\langle y^*, \bar{x}_M\rangle}{2tr}$$

$$= \frac{\langle y^*_M, \bar{x}_M\rangle}{2r}.$$

This implies

$$\theta^* \in \widehat{D}^* P_{r\mathbb{C}_M}(\bar{x})(y^*) \implies \langle y^*_M, \bar{x}_M\rangle \leq 0. \tag{4.18}$$

By (4.15.2), we must have $-(J^*(y^*))_M \in \bar{x}_r^{\Uparrow}$. Then, there is $\delta > 0$ such that $\|\bar{x} - t(J^*(y^*))_M\|_p > r$, for all $t \in (0, \delta)$. In this case, we take a directional line segment in the limit (4.16) as $u = \bar{x} - t(J^*(y^*))_M$, for $t \downarrow 0$ with $t < \delta$. By (4.5), we have

$$P_{r\mathbb{C}_M}(\bar{x} - t(J^*(y^*))_M) = \frac{r}{\|(\bar{x} - t(J^*(y^*))_M)_M\|_p}(\bar{x} - t(J^*(y^*))_M)_M + (\bar{x} - t(J^*(y^*))_M)_{\overline{M}}$$

$$= \frac{r}{\|\bar{x}_M - t(J^*(y^*))_M\|_p}(\bar{x}_M - t(J^*(y^*))_M) + \bar{x}_{\overline{M}}.$$

We have

$$0 \geq \limsup_{u \to \bar{x}} \frac{-\langle y^*, P_{r\mathbb{C}_M}(u) - \bar{x}\rangle}{\|u - \bar{x}\|_p + \|P_{r\mathbb{C}_M}(u) - \bar{x}\|_p}$$

$$\geq \limsup_{t \downarrow 0, t < \delta} \frac{-\langle y^*, \frac{r}{\|\bar{x}_M - t(J^*(y^*))_M\|_p}(\bar{x}_M - t(J^*(y^*))_M) + \bar{x}_{\overline{M}} - \bar{x}\rangle}{\left\|\bar{x} - t(J^*(y^*))_M - \bar{x}\right\|_p + \left\|\frac{r}{\|\bar{x}_M - t(J^*(y^*))_M\|_p}(\bar{x}_M - t(J^*(y^*))_M) + \bar{x}_{\overline{M}} - \bar{x}\right\|_p}$$

$$= \limsup_{t \downarrow 0, t < \delta} \frac{-\langle y^*, \frac{r}{\|\bar{x}_M - t(J^*(y^*))_M\|_p}(\bar{x}_M - t(J^*(y^*))_M) - \bar{x}_M\rangle}{t\|(J^*(y^*))_M\|_p + \left\|\frac{r}{\|\bar{x}_M - t(J^*(y^*))_M\|_p}(\bar{x}_M - t(J^*(y^*))_M) - \bar{x}_M\right\|_p}$$

$$= \limsup_{t \downarrow 0, t < \delta} \frac{-\langle y^*, \left(\frac{r}{\|\bar{x}_M - t(J^*(y^*))_M\|_p} - 1\right)\bar{x}_M - \frac{tr(J^*(y^*))_M}{\|\bar{x}_M - t(J^*(y^*))_M\|_p}\rangle}{t\|(J^*(y^*))_M\|_p + \left\|\left(\frac{r}{\|\bar{x}_M - t(J^*(y^*))_M\|_p} - 1\right)\bar{x}_M - \frac{tr(J^*(y^*))_M}{\|\bar{x}_M - t(J^*(y^*))_M\|_p}\right\|_p}$$

$$= \operatorname*{limsup}_{t\downarrow 0, t<\delta} \frac{\frac{tr\langle y^*, (J^*(y^*))_M \rangle}{\|\bar{x}_M - t(J^*(y^*))_M\|_p} - \langle y^*, \left(\frac{r}{\|\bar{x}_M - t(J^*(y^*))_M\|_p} - 1\right)\bar{x}_M\rangle}{t\|(J^*(y^*))_M\|_p + \left\|\left(\frac{r}{\|\bar{x}_M - t(J^*(y^*))_M\|_p} - 1\right)\bar{x}_M - \frac{tr(J^*(y^*))_M}{\|\bar{x}_M - t(J^*(y^*))_M\|_p}\right\|_p}$$

$$= \operatorname*{limsup}_{t\downarrow 0, t<\delta} \frac{\frac{tr\langle y^*, (J^*(y^*))_M \rangle}{\|\bar{x}_M - t(J^*(y^*))_M\|_p} + \frac{tr\langle y^*, \bar{x}_M\rangle}{\|\bar{x}_M - t(J^*(y^*))_M\|_p \|\bar{x}_M\|_p}\left(\frac{\|\bar{x}_M - t(J^*(y^*))_M\|_p - \|\bar{x}_M\|_p}{t}\right)}{t\|(J^*(y^*))_M\|_p + \left\|\frac{tr}{\|\bar{x}_M - t(J^*(y^*))_M\|_p \|\bar{x}_M\|_p}\left(\frac{\|\bar{x}_M - t(J^*(y^*))_M\|_p - \|\bar{x}_M\|_p}{t}\right)\bar{x}_M + \frac{tr(J^*(y^*))_M}{\|\bar{x}_M - t(J^*(y^*))_M\|_p}\right\|_p}$$

$$= \operatorname*{limsup}_{t\downarrow 0, t<\delta} \frac{\frac{tr\langle y_M^*, (J^*(y^*))_M \rangle}{\|\bar{x}_M - t(J^*(y^*))_M\|_p} + \frac{t\langle y^*, \bar{x}_M\rangle}{\|\bar{x}_M - t(J^*(y^*))_M\|_p}\left(\frac{\|\bar{x}_M - t(J^*(y^*))_M\|_p - \|\bar{x}_M\|_p}{t}\right)}{t\|(J^*(y^*))_M\|_p + \frac{t}{\|\bar{x}_M - t(J^*(y^*))_M\|_p}\left\|\left(\frac{\|\bar{x}_M - t(J^*(y^*))_M\|_p - \|\bar{x}_M\|_p}{t}\right)\bar{x}_M + (J^*(y^*))_M\right\|_p}$$

$$= \operatorname*{limsup}_{t\downarrow 0, t<\delta} \frac{\frac{r\|y_M^*\|_q^2}{\|\bar{x}_M - t(J^*(y^*))_M\|_p} + \frac{\langle y^*, \bar{x}_M\rangle}{\|\bar{x}_M - t(J^*(y^*))_M\|_p}\left(\frac{\|\bar{x}_M - t(J^*(y^*))_M\|_p - \|\bar{x}_M\|_p}{t}\right)}{\|(J^*(y^*))_M\|_p + \frac{1}{\|\bar{x}_M - t(J^*(y^*))_M\|_p}\left\|\left(\frac{\|\bar{x}_M - t(J^*(y^*))_M\|_p - \|\bar{x}_M\|_p}{t}\right)\bar{x}_M + (J^*(y^*))_M\right\|_p}$$

$$= \frac{\frac{r\|y_M^*\|_q^2}{\|\bar{x}_M\|_p} + \frac{\langle y_M^*, \bar{x}_M\rangle}{\|\bar{x}_M\|_p} \frac{\langle J(\bar{x}_M), (-J^*(y^*))_M\rangle}{\|\bar{x}_M\|_p}}{\|(J^*(y^*))_M\|_p + \frac{1}{\|\bar{x}_M\|_p}\left\|\frac{\langle J(\bar{x}_M), (-J^*(y^*))_M\rangle}{\|\bar{x}_M\|_p}\bar{x}_M + (J^*(y^*))_M\right\|_p}$$

$$= \frac{\|y_M^*\|_q^2 - \frac{\langle y_M^*, \bar{x}_M\rangle\langle J(\bar{x}_M), (J^*(y^*))_M\rangle}{r^2}}{\|(J^*(y^*))_M\|_p + \frac{1}{\|\bar{x}_M\|_p}\left\|\frac{\langle J(\bar{x}_M), (-J^*(y^*))_M\rangle}{\|\bar{x}_M\|_p}\bar{x}_M + (J^*(y^*))_M\right\|_p}$$

$\geq 0.$

This implies

$$\|y_M^*\|_q^2 - \frac{\langle y_M^*, \bar{x}_M\rangle\langle J(\bar{x}_M), (J^*(y^*))_M\rangle}{r^2} = 0. \tag{4.19}$$

Since $\|J(\bar{x}_M)\|_q = \|\bar{x}_M\|_p = r$, by (4.18), to satisfy (4.19), we must have

$$\langle y_M^*, \bar{x}_M\rangle = -\|y_M^*\|_q \, r = -\|y_M^*\|_q \|\bar{x}_M\|_p, \tag{4.20}$$

and $\langle J(\bar{x}_M), (J^*(y^*))_M\rangle = -r \left\|(J^*(y^*))_M\right\|_p = -\|J(\bar{x}_M)\|_q \left\|(J^*(y^*))_M\right\|_p = -\|y_M^*\|_q \|\bar{x}_M\|_p$.

(4.20) proves (4.15.3). Hence (4.15) is proved. Then, we prove the reverse of (4.15); that is

$$\langle y_M^*, \bar{x}_M \rangle = -\|y_M^*\|_q \|\bar{x}_M\|_p, \; y_{\overline{M}}^* = \theta^* \text{ and } -(J^*(y^*))_M \notin \bar{x}_r^{\Downarrow}$$
$$\Rightarrow \theta^* \in \widehat{D}^* P_{r\mathbb{B}}(\bar{x})(y^*). \tag{4.21}$$

By $\|\bar{x}_M\|_p = r$ and $\left\|\left(P_{r\mathbb{C}_M}(u)\right)_M\right\| \leq r$, for all $u \in l_p$, we calculate

$$\limsup_{u \to \bar{x}} \frac{-\langle y^*, \; P_{r\mathbb{C}_M}(u) - \bar{x} \rangle}{\|u - \bar{x}\|_p + \|P_{r\mathbb{C}_M}(u) - \bar{x}\|_p}$$

$$= \limsup_{u \to \bar{x}} \frac{\langle y^*, \bar{x} \rangle - \langle y^*, \; P_{r\mathbb{C}_M}(u) \rangle}{\|u - \bar{x}\|_p + \|P_{r\mathbb{C}_M}(u) - \bar{x}\|_p}$$

$$= \limsup_{u \to \bar{x}} \frac{\langle y_M^*, \bar{x} \rangle - \langle y_M^*, \; P_{r\mathbb{C}_M}(u) \rangle}{\|u - \bar{x}\|_p + \|P_{r\mathbb{C}_M}(u) - \bar{x}\|_p}$$

$$= \limsup_{u \to \bar{x}} \frac{\langle y_M^*, \bar{x}_M \rangle - \langle y_M^*, \; P_{r\mathbb{C}_M}(u) \rangle}{\|u - \bar{x}\|_p + \|P_{r\mathbb{C}_M}(u) - \bar{x}\|_p}$$

$$= \limsup_{u \to \bar{x}} \frac{-\|y_M^*\|_q \|\bar{x}_M\|_p - \langle y_M^*, \; (P_{r\mathbb{C}_M}(u))_M \rangle}{\|u - \bar{x}\|_p + \|P_{r\mathbb{C}_M}(u) - \bar{x}\|_p}$$

$$= \limsup_{u \to \bar{x}} \frac{-\|y_M^*\|_q \, r - \langle y_M^*, \; (P_{r\mathbb{C}_M}(u))_M \rangle}{\|u - \bar{x}\|_p + \|P_{r\mathbb{C}_M}(u) - \bar{x}\|_p}$$

$$\leq \limsup_{u \to \bar{x}} \frac{-\|y_M^*\|_q \, r + \|y_M^*\|_q \left\|(P_{r\mathbb{C}_M}(u))_M\right\|_p}{\|u - \bar{x}\|_p + \|P_{r\mathbb{C}_M}(u) - \bar{x}\|_p}$$

$$\leq \limsup_{u \to \bar{x}} \frac{\|y_M^*\|_q \left(\left\|(P_{r\mathbb{C}_M}(u))_M\right\|_p - r\right)}{\|u - \bar{x}\|_p + \|P_{r\mathbb{C}_M}(u) - \bar{x}\|_p}$$

$$\leq 0.$$

By (4.16), this implies $\theta^* \in \widehat{D}^* P_{r\mathbb{C}_M}(\bar{x})(y^*)$, which proves (4.21). By (4.15) and (4.21), part (b) in (iii) is proved.

Proof of (c). The proof of part (c) in (iii) in this theorem is similar to the proof of part (c) in (iii) of Theorem 3.1. Let $\bar{x}^* = J(\bar{x})$ with $\|\bar{x}_M\|_p = r$. Let $z^* \in l_q$, we have

$$z^* \in \widehat{D}^* P_{r\mathbb{C}_M}(\bar{x})(\bar{x}^*) \iff \limsup_{u \to \bar{x}} \frac{\langle (z^*, -\bar{x}^*), \; (u, P_{r\mathbb{C}_M}(u) - (\bar{x}, P_{r\mathbb{C}_M}(\bar{x})) \rangle}{\|u - \bar{x}\|_p + \|P_{r\mathbb{C}_M}(u) - P_{r\mathbb{C}_M}(\bar{x})\|_p} \leq 0$$

$$\iff \limsup_{u \to \bar{x}} \frac{\langle (z^*, -\bar{x}^*), \; (u, P_{r\mathbb{C}_M}(u) - (\bar{x}, \bar{x})) \rangle}{\|u - \bar{x}\|_p + \|P_{r\mathbb{C}_M}(u) - \bar{x}\|_p} \leq 0$$

$$\iff \limsup_{u \to \bar{x}} \frac{\langle z^*, \; u - \bar{x} \rangle - \langle \bar{x}^*, \; P_{r\mathbb{C}_M}(u) - \bar{x} \rangle}{\|u - \bar{x}\|_p + \|P_{r\mathbb{C}_M}(u) - \bar{x}\|_p} \leq 0. \tag{4.22}$$

We take a directional line segment in the limit in (4.22) $u = \bar{x} - t\bar{x}_M$, for $t \downarrow 0$ with $t < 1$. Then,

$$\|(\bar{x} - t\bar{x}_M)_M\|_p = \|\bar{x}_M - t\bar{x}_M\|_p = (1-t)\|\bar{x}_M\|_p < r, \text{ for all } 0 < t < 1.$$

By (4.5), this implies

$$P_{r\mathbb{C}_M}(\bar{x} - t\bar{x}_M) = \bar{x} - t\bar{x}_M, \text{ for all } 0 < t < 1.$$

By (4.22) and $\langle z_{\bar{M}}^*, \bar{x}_M \rangle = 0$, we have

$$0 \geq \limsup_{u \to \bar{x}} \frac{\langle z^*, u-\bar{x}\rangle - \langle \bar{x}^*, P_{r\mathbb{C}_M}(u)-\bar{x}\rangle}{\|u-\bar{x}\|_p + \|P_{r\mathbb{C}_M}(u)-\bar{x}\|_p}$$

$$\geq \limsup_{t \downarrow 0, t<1} \frac{\langle z^*, \bar{x}-t\bar{x}_M-\bar{x}\rangle - \langle \bar{x}^*, \bar{x}-t\bar{x}_M-\bar{x}\rangle}{\|\bar{x}-t\bar{x}_M-\bar{x}\|_p + \|\bar{x}-t\bar{x}_M-\bar{x}\|}$$

$$= \limsup_{t \downarrow 0, t<1} \frac{\langle z_M^* + z_{\bar{M}}^*, -t\bar{x}_M\rangle - \langle \bar{x}^*, -t\bar{x}_M\rangle}{2t\|\bar{x}\|_p}$$

$$= \frac{-\langle z_M^*, t\bar{x}_M\rangle + \|\bar{x}_M\|_p^2}{2\|\bar{x}_M\|_p}$$

$$= \frac{-\langle z_M^*, t\bar{x}_M\rangle + r^2}{2r}.$$

This concludes that

$$z^* \in \widehat{D}^* P_{r\mathbb{C}_M}(\bar{x})(J(\bar{x})) \implies \langle z_M^*, t\bar{x}_M\rangle \geq r^2. \tag{4.23}$$

On the other hand, we take a different directional line segment in the limit in (4.22) as $u = \bar{x} + t\bar{x}_M$, for $t \downarrow 0$ with $t < 1$. Then,

$$\|(\bar{x} + t\bar{x}_M)_M\|_p = \|\bar{x}_M + t\bar{x}_M\|_p = (1+t)\|\bar{x}_M\|_p > r, \text{ for all } 0 < t < 1.$$

By (4.5) and $\|\bar{x}_M\|_p = r$, this implies

$$P_{r\mathbb{C}_M}(\bar{x} - t\bar{x}_M) = \frac{r}{\|(\bar{x}+t\bar{x}_M)_M\|_p}(\bar{x}+t\bar{x}_M)_M + (\bar{x}+t\bar{x}_M)_{\bar{M}}$$

$$= \frac{r}{\|\bar{x}_M + t\bar{x}_M\|_p}(\bar{x}_M + t\bar{x}_M) + \bar{x}_{\bar{M}}$$

$$= \bar{x}_M + \bar{x}_{\bar{M}}$$

$$= \bar{x}, \text{ for all } 0 < t < 1.$$

By (4.22), we have

$$0 \geq \limsup_{u \to \bar{x}} \frac{\langle z^*, u-\bar{x}\rangle - \langle \bar{x}^*, P_{r\mathbb{C}_M}(u)-\bar{x}\rangle}{\|u-\bar{x}\|_p + \|P_{r\mathbb{C}_M}(u)-\bar{x}\|_p}$$

$$\geq \limsup_{t\downarrow 0} \frac{\langle z^*, \ \bar{x}+t\bar{x}_M-\bar{x}\rangle - \langle \bar{x}^*, \ \bar{x}-\bar{x}\rangle}{\|\bar{x}+t\bar{x}_M-\bar{x}\|_p + \|\bar{x}-\bar{x}\|_p}$$

$$= \limsup_{t\downarrow 0} \frac{\langle z_M^*+z_{\bar{M}}^*, \ t\bar{x}_M\rangle - \langle \bar{x}^*,\theta\rangle}{t\|\bar{x}_M\|_p}$$

$$= \frac{\langle z_M^*, \ \bar{x}_M\rangle}{\|\bar{x}\|_p}$$

$$= \frac{\langle z_M^*, \ \bar{x}_M\rangle}{r}.$$

This implies that

$$z^* \in \widehat{D}^* P_{r\mathbb{C}_M}(\bar{x})(J(\bar{x})) \implies \langle z_M^*, \ \bar{x}_M\rangle \leq 0. \tag{4.24}$$

(4.24) contradicts to (4.23). This implies that $\widehat{D}^* P_{r\mathbb{C}_M}(\bar{x})(J(\bar{x})) = \emptyset$, which proves (c). $\square$

In particular, when $M = \mathbb{N}$, which implies $\bar{M} = \emptyset$, then

$$x_M = x, \ x_{\bar{M}} = \theta, \text{ for any } x \in l_p \quad \text{and} \quad y_M^* = y^*, y_{\bar{M}}^* = \theta^*, \text{ for any } y^* \in l_q.$$

The convex and closed cylinder $r\mathbb{C}_{\mathbb{N}}$ becomes the closed ball $r\mathbb{B}$. That is, $r\mathbb{C}_{\mathbb{N}} = r\mathbb{B}$.

As a consequence of Theorem 4.1, we obtain the major parts (i), (ii), (a, c) of (iii) in Theorem 3.1 immediately. We get some new result, which is similar to part (b) of (iii) in Theorem 3.1.

**Corollary 4.2.** *For any $r > 0$, Mordukhovich derivatives of the metric projection $P_{r\mathbb{B}}$ at an arbitrarily given point $\bar{x} \in r\mathbb{S}$ satisfies that, for $y^* \in l_q\setminus\{\theta\}$,*

$$\theta^* \in \widehat{D}^* P_{r\mathbb{B}}(\bar{x})(y^*) \iff -J^*(y^*) \notin \bar{x}_r^{\Downarrow} \text{ and } \langle y^*,\bar{x}\rangle = -r\|y^*\|_q.$$

*In particular,*

$$\theta^* \in \widehat{D}^* P_{r\mathbb{C}_M}(\bar{x})(-J(\bar{x})).$$

5. **Mordukhovich derivatives of the metric projection onto the positive cone in real Banach space $L_p(S)$**

Let $(S, \mathcal{A}, \mu)$ be a positive and complete measure space. In this section, we consider the real uniformly convex and uniformly smooth Banach space $(L_p(S), \|\cdot\|_p)$ with dual space $(L_q(S), \|\cdot\|_q)$, in which $p$ and $q$ satisfy $1 < p, q < \infty$ and $\frac{1}{p} + \frac{1}{q} = 1$. $L_p(S)$ and $L_q(S)$ share the same origin $\theta = \theta^*$. For the sake of distinction between $L_p(S)$ and its dual space $L_q(S)$, we use English letters $f, g, h, \ldots$ for the elements in $L_p(S)$, and we use Greek letters $\varphi, \psi, \xi, \ldots$ for the elements in the dual space $L_q(S)$. We define the positive cone $K_p$ in $L_p(S)$ and investigate the Mordukhovich derivatives of the metric projection onto $K_p$.

The normalized duality mapping $J: L_p(S) \to L_q(S)$ has the following representations, for any given $f \in L_p(S)$ with $f \neq \theta$,

$$(Jf)(s) = \frac{|f(s)|^{p-1}\operatorname{sign}(f(s))}{\|f\|_p^{p-2}} = \frac{|f(s)|^{p-2}f(s)}{\|f\|_p^{p-2}}, \quad \text{for all } s \in S. \tag{5.1}$$

We define a subset $K_p$ of $L_p(S)$ as follows:

$$K_p = \{f \in L_p(S) : f(s) \geq 0, \text{ for } \mu\text{-almost all } s \in S\}.$$

Similarly, in $L_q(S)$, which is the dual space of $L_p(S)$, we define

$$K_q = \{\varphi \in L_q(S) : \varphi(s) \geq 0, \text{ for } \mu\text{-almost all } s \in S\}.$$

$K_p$ and $K_q$ are pointed closed and convex cones, which are called the positive cones of $L_p(S)$ and $L_q(S)$, respectively. Let $\preccurlyeq_p$ and $\preccurlyeq_q$ be the two partial orders on $L_p(S)$ and $L_q(S)$ induced by the two pointed closed and convex cones $K_p$ and $K_q$, respectively. We define the ordered intervals in $L_p(S)$ and $L_q(S)$. For any $f, g \in L_p(S)$ with $f \preccurlyeq_p g$, and for any $\xi, \psi \in L_q(S)$ with $\varphi \preccurlyeq_p \psi$, we write

$$[f, g]_{\preccurlyeq_p} = \{h \in L_p(S) : f \preccurlyeq_p h \preccurlyeq_p g\},$$

$$[\xi, \psi]_{\preccurlyeq_q} = \{\varphi \in L_q(S) : \xi \preccurlyeq_q \varphi \preccurlyeq_q \psi\}.$$

**Lemma 5.1 in [15]**. *Both $K_p$ and $K_q$ have empty interior.*

**Lemma 5.2 in [15]**. *The metric projection $P_{K_p} : L_p(S) \to K_p$ has the following representations.*

$$(P_{K_p}f)(s) = \begin{cases} f(s), & \text{if } f(s) > 0, \\ 0, & \text{if } f(s) \leq 0, \end{cases} \quad \text{for any } f \in L_p(S). \tag{5.2}$$

**Lemma 5.3 in [15]**. *The metric projection $P_{K_p} : L_p(S) \to K_p$ has the following properties*

(a) $P_{K_p}(f) = f$, *for any $f \in K_p$;*
(b) $P_{K_p}(f) = \theta$, *for any $f \in -K_p$;*
(c) $P_{K_p}(f + g) = f + g$, *for any $f, g \in K_p$;*
(d) *For any $f \in L_p(S)$,*

$$P_{K_p}(\lambda f) = \lambda P_{K_p}(f), \text{ for any } \lambda \geq 0. \tag{5.3}$$

**Theorem 5.4 in [15]**. *$P_{K_p}$ is not Fréchet differentiable at every point in $L_p(S)$, that is,*

$$\nabla P_{K_p}(f) \text{ does not exist, for any } f \in L_p(S).$$

We need the following notations. For any $f \in L_p(S)$, for all $s \in S$, we write

$$f^+(s) = \begin{cases} f(s), & \text{for } f(s) > 0, \\ 0, & \text{for } f(s) \leq 0. \end{cases}$$

and

$$f^-(s) = \begin{cases} f(s), & \text{for } f(s) < 0, \\ 0, & \text{for } f(s) \geq 0. \end{cases}$$

**Lemma 5.1.** *For any $f \in L_p(S) \setminus \{\theta\}$, $f^+$ and $f^-$ have the following properties.*

(a) $f = f^+ + f^-$;
(b) $P_{K_p}(f) = f^+$;
(c) $f = P_{K_p}(f) + f^-$;
(d) $(J(f))^+ = \dfrac{\|f^+\|_p^{p-2}}{\|f\|_p^{p-2}} J(f^+)$;
(e) $(J(f))^- = \dfrac{\|f^-\|_p^{p-2}}{\|f\|_p^{p-2}} J(f^-)$;
(f) $\langle J(f), f^+ \rangle = \langle (J(f))^+, f^+ \rangle = \dfrac{\|f^+\|_p^{p-2}}{\|f\|_p^{p-2}} \langle J(f^+), f^+ \rangle = \dfrac{\|f^+\|_p^p}{\|f\|_p^{p-2}}$;
(g) $\langle J(f), f^- \rangle = \langle (J(f))^-, f^- \rangle = \dfrac{\|f^-\|_p^{p-2}}{\|f\|_p^{p-2}} \langle J(f^-), f^- \rangle = \dfrac{\|f^-\|_p^p}{\|f\|_p^{p-2}}$.

*Proof.* We only prove part (d). Rest of the proof of this lemma can be similarly proved. If $\|f^+\|_p = 0$, then, (d) is clear. So, we suppose $\|f^+\|_p \neq 0$. By (5.1), notice that,

$$J(f)(s) > 0 \quad \Leftrightarrow \quad f(s) > 0, \text{ for every } s \in S.$$

Then, by (5.1) again, we have

$$(J(f))^+(s) = \begin{cases} J(f)(s), & \text{if } J(f)(s) > 0 \\ 0, & \text{if } J(f)(s) \leq 0 \end{cases}$$

$$= \begin{cases} J(f)(s), & \text{if } f(s) > 0 \\ 0, & \text{if } f(s) \leq 0 \end{cases}$$

$$= \begin{cases} \dfrac{|f(s)|^{p-1} \operatorname{sign}(f(s))}{\|f\|_p^{p-2}}, & \text{if } f(s) > 0 \\ 0, & \text{if } f(s) \leq 0 \end{cases}$$

$$= \begin{cases} \dfrac{|f(s)|^{p-1}}{\|f\|_p^{p-2}}, & \text{if } f(s) > 0 \\ 0, & \text{if } f(s) \leq 0 \end{cases}$$

$$= \dfrac{\|f^+\|_p^{p-2}}{\|f\|_p^{p-2}} \begin{cases} \dfrac{|f(s)|^{p-1}}{\|f^+\|_p^{p-2}}, & \text{if } f(s) > 0 \\ 0, & \text{if } f(s) \leq 0 \end{cases}$$

$$= \frac{\|f^+\|_p^{p-2}}{\|f\|_p^{p-2}} \begin{cases} \frac{|f^+(s)|^{p-1}}{\|f^+\|_p^{p-2}}, & \text{if } f(s) > 0 \\ 0, & \text{if } f(s) \le 0 \end{cases}$$

$$= \frac{\|f^+\|_p^{p-2}}{\|f\|_p^{p-2}} (J(f^+))(s), \text{ for all } s \in S. \qquad \square$$

**Theorem 5.2.** *The metric projection* $P_{K_p}: L_p(S) \to K_p$ *has the following Mordukhovich derivatives. For any* $f \in L_p(S)$,

(a) $\widehat{D}^* P_{K_p}(f)(\theta^*) = \{\theta^*\}$;

(b) *For any* $\varphi \in l_q$,

$$\theta^* \in \widehat{D}^* P_{K_p}(f)(\varphi)$$

$$\iff \mu(\{s \in S : \varphi(s) \neq 0 \text{ and } f(s) > 0\} \cup \{s \in S : \varphi(s) < 0 \text{ and } f(s) \le 0\}) = 0.$$

*It implies that*,

(b₁) *For any* $f \in -K_p$, *we have*

$$\theta^* \in \widehat{D}^* P_{K_p}(f)(\varphi), \text{ for any } \varphi \in K_q;$$

(b₂) *For any* $f \in K_p \setminus \{\theta\}$,

$$\theta^* \notin \widehat{D}^* P_{K_p}(f)(J(f));$$

(c) $J(f) \in \widehat{D}^* P_{K_p}(f)(J(f))$, *for any* $f \in K_p$;

(d) *For any given* $\psi \in K_q$, *we have*

$$\widehat{D}^* P_{K_p}(\theta)(\psi) = [\theta^*, \psi]_{\leqslant_q}.$$

*Proof.* Proof of (a). For any $f \in L_p(S)$, it is clear that

$$\theta^* \in \widehat{D}^* P_{K_p}(f)(\theta^*). \tag{5.8}$$

Next, we show $\varphi \notin \widehat{D}^* P_{K_p}(f)(\theta^*)$, for any $\varphi \in L_q(S) \setminus \{\theta^*\}$. By (2.20), we have

$$\varphi \in \widehat{D}^* P_{K_p}(f)(\theta^*) \iff \limsup_{h \to f} \frac{\langle (\varphi, -\theta^*), (h, P_{K_p}(h)) - (f, P_{K_p}(f)) \rangle}{\|h-f\|_p + \|P_{K_p}(h) - P_{K_p}(f)\|_p} \le 0$$

$$\iff \limsup_{h \to f} \frac{\langle \varphi, h-f \rangle}{\|h-f\|_p + \|P_{K_p}(h) - P_{K_p}(f)\|_p} \le 0. \tag{5.9}$$

Since $\varphi \in L_q(S)\setminus\{\theta^*\}$, then, we have either $\varphi^+ \neq \theta^*$ or $\varphi^- \neq \theta^*$ or both $\varphi^+ \neq \theta^*$ and $\varphi^- \neq \theta^*$. We consider it case by case.

Case 1. $\varphi^+ \neq \theta^*$. This implies that $(J^*(\varphi))^+ \neq \theta$. In this case, we take a directional line segment in the limit (5.9) as $h = f + t(J^*(\varphi))^+$, for $t \downarrow 0$ with $t < 1$. One can check that,

$$\left(f(s) + t(J^*(\varphi))^+(s)\right)^+ \geq f(s) + t(J^*(\varphi))^+(s), \text{ for all } s \in S.$$

and $\quad f^+(s) + t(J^*(\varphi))^+(s) \geq \left(f(s) + t(J^*(\varphi))^+(s)\right)^+(s) \geq f^+(s),$ for all $s \in S$.

This reduces to

$$0 \leq \left(f(s) + t(J^*(\varphi))^+(s)\right)^+(s) - f^+(s) \leq f^+(s) + t(J^*(\varphi))^+(s) - f^+(s), \text{ for all } s \in S.$$

By Lemma 5.1, this implies

$$0 \geq \limsup_{h \to f} \frac{\langle \varphi, h-f \rangle}{\|h-f\|_p + \|P_{K_p}(h) - P_{K_p}(f)\|_p}$$

$$\geq \limsup_{t \downarrow 0} \frac{\langle \varphi, f+t(J^*(\varphi))^+ - f \rangle}{\|f+t(J^*(\varphi))^+ - f\|_p + \|(f+t(J^*(\varphi))^+)^+ - f^+\|_p}$$

$$= \limsup_{t \downarrow 0} \frac{\langle \varphi, t(J^*(\varphi))^+ \rangle}{\|f+t(J^*(\varphi))^+ - f\|_p + \|(f+t(J^*(\varphi))^+)^+ - f^+\|_p}$$

$$= \limsup_{t \downarrow 0} \frac{t\langle \varphi, \frac{\|\varphi^+\|_q^{q-2}}{\|\varphi\|_q^{q-2}} J^*(\varphi^+)\rangle}{\|f+t(J^*(\varphi))^+ - f\|_p + \|(f+t(J^*(\varphi))^+)^+ - f^+\|_p}$$

$$= \limsup_{t \downarrow 0} \frac{t\frac{\|\varphi^+\|_q^{q-2}}{\|\varphi\|_q^{q-2}}\langle \varphi^+, J^*(\varphi^+)\rangle}{\|f+t(J^*(\varphi))^+ - f\|_p + \|(f+t(J^*(\varphi))^+)^+ - f^+\|_p}$$

$$= \limsup_{t \downarrow 0} \frac{t\frac{\|\varphi^+\|_q^q}{\|\varphi\|_q^{q-2}}}{\|f+t(J^*(\varphi))^+ - f\|_p + \left\|\left(f+t(J^*(\varphi))^+\right)^+ - f^+\right\|_p}$$

$$\geq \limsup_{t \downarrow 0} \frac{t\frac{\|\varphi^+\|_q^q}{\|\varphi\|_q^{q-2}}}{t\|(J^*(\varphi))^+\|_p + \|f^+ + t(J^*(\varphi))^+ - f^+\|_p}$$

$$= \limsup_{t \downarrow 0} \frac{t\frac{\|\varphi^+\|_q^q}{\|\varphi\|_q^{q-2}}}{2t\|(J^*(\varphi))^+\|_p}$$

$$= \limsup_{t\downarrow 0} \frac{t\frac{\|\varphi^+\|_q^q}{\|\varphi\|_q^{q-2}}}{2t\frac{\|\varphi^+\|_q^{q-2}}{\|\varphi\|_q^{q-2}}\|J^*(\varphi^+)\|_p}$$

$$= \frac{\frac{\|\varphi^+\|_q^q}{\|\varphi\|_q^{q-2}}}{2\frac{\|\varphi^+\|_q^{q-2}}{\|\varphi\|_q^{q-2}}\|\varphi^+\|_q}$$

$$= \frac{\|\varphi^+\|_q}{2}$$

$$> 0. \tag{5.10}$$

(5.10) is a contradiction.

Case 2. $\varphi^- \neq \theta^*$. This implies that $(J^*(\varphi))^- \neq \theta$. In this case, we take a directional line segment in the limit (5.9) as, $h = f + t(J^*(\varphi))^-$, for $t \downarrow 0$ with $t < 1$. One can check that

$$f^+(s) + t(J^*(\varphi))^-(s) \leq (f + t(J^*(\varphi))^-)^+(s) \leq f^+(s), \text{ for all } s \in S.$$

This is,

$$0 \leq f^+(s) - (f + t(J^*(\varphi))^-)^+(s) \leq f^+(s) - (f^+(s) + t(J^*(\varphi))^-(s)), \text{ for all } s \in S.$$

By Lemma 5.1, similar to case 1, this implies

$$0 \geq \limsup_{h \to f} \frac{\langle \varphi, P_{K_p}(h)-f \rangle}{\|h-f\|_p + \|P_{K_p}(h)-P_{K_p}(f)\|_p}$$

$$\geq \limsup_{t \downarrow 0} \frac{\langle \varphi, f+t(J^*(\varphi))^- - f \rangle}{\|f+t(J^*(\varphi))^- - f\|_p + \|P_{K_p}(f+t(J^*(\varphi))^-)-P_{K_p}(f)\|_p}$$

$$= \limsup_{t \downarrow 0} \frac{t\langle \varphi, (J^*(\varphi))^- \rangle}{t\|(J^*(\varphi))^-\|_p + \|(f+t(J(g))^-)^+ - f^+\|_p}$$

$$= \limsup_{t \downarrow 0} \frac{t\frac{\|\varphi^-\|_q^q}{\|\varphi\|_q^{q-2}}}{t\|(J^*(\varphi))^-\|_p + \|f^+ - (f+t(J^*(\varphi))^-)^+\|_p}$$

$$\geq \limsup_{t \downarrow 0} \frac{t\frac{\|\varphi^-\|_q^q}{\|\varphi\|_q^{q-2}}}{t\|(J^*(\varphi))^-\|_p + \|f^+ - (f^+ + t(J^*(\varphi))^-)\|_p}$$

$$\geq \limsup_{t \downarrow 0} \frac{t\frac{\|\varphi^-\|_q^q}{\|\varphi\|_q^{q-2}}}{2t\|(J^*(\varphi))^-\|_p}$$

$$= \frac{\|\varphi^-\|_q}{2}$$

$$> 0. \tag{5.11}$$

This is a contradiction. By (5.10) and (5.11), it proved

$$\varphi \in L_q(S)\setminus\{\theta^*\} \implies \varphi \notin \widehat{D}^* P_{K_p}(f)(\theta^*). \tag{5.12}$$

Part (a) is proved by (5.8) and (5.12).

Proof of (b). We firstly prove the part "$\Rightarrow$". The proof of this part "$\Rightarrow$" in (b) is divided to the following three parts.

$$\theta^* \in \widehat{D}^* P_{K_p}(f)(\varphi)$$

$$\Rightarrow \quad \begin{array}{l} (B_1) \ \mu(\{s \in S: \varphi(s) < 0 \text{ and } f(s) > 0\}) = 0, \\ (B_2) \ \mu(\{s \in S: \varphi(s) > 0 \text{ and } f(s) > 0\}) = 0, \\ \text{and } (B_3) \ \mu(\{s \in S: \varphi(s) < 0 \text{ and } f(s) \leq 0\}) = 0. \end{array}$$

Let $f \in L_p(S)$. For $\varphi \in L_q(S)\setminus\{\theta^*\}$, by (2.20) and Lemma 5.1, we have

$$\theta^* \in \widehat{D}^* P_{K_p}(f)(\varphi) \iff \limsup_{h \to f} \frac{\langle(\theta^*, -\varphi), (h, P_{K_p}(h)) - (f, P_{K_p}(f))\rangle}{\|h-f\|_p + \|P_{K_p}(h) - P_{K_p}(f)\|_p} \leq 0$$

$$\iff \limsup_{h \to f} \frac{-\langle\varphi, P_{K_p}(h) - P_{K_p}(f)\rangle}{\|h-f\|_p + \|P_{K_p}(h) - P_{K_p}(f)\|_p} \leq 0$$

$$\iff \limsup_{h \to f} \frac{-\langle\varphi, h^+ - f^+\rangle}{\|h-f\|_p + \|h^+ - f^+\|_p} \leq 0. \tag{5.13}$$

Proof of (B$_1$). Assume, by contradiction that $\mu(\{s \in S: \varphi(s) < 0\} \cap \{s \in S: f(s) > 0\}) > 0$. By the properties of $J^*$, we have that $\varphi(s) < 0$ if and only if $J^*(\varphi)(s) < 0$. This implies

$$\mu(\{s \in S: J^*(\varphi)(s) < 0\} \cap \{s \in S: f(s) > 0\}) > 0.$$

In this case, let $A = \{s \in S: J^*(\varphi)(s) < 0\} \cap \{s \in S: f(s) > 0\}$. Define

$$\varphi_A(s) = \begin{cases} \varphi(s), & \text{if } s \in A, \\ 0, & \text{if } s \notin A. \end{cases}$$

We take a directional line segment in the limit (5.13) as $h = f - t(J^*(\varphi_A))$, for $t \downarrow 0$. One has

$$(f - t(J^*(\varphi_A)))^+ = f^+ - t(J^*(\varphi_A)).$$

This implies

$$0 \geq \limsup_{h \to f} \frac{-\langle\varphi, h^+ - f^+\rangle}{\|h-f\|_p + \|h^+ - f^+\|_p}$$

$$\geq \limsup_{t\downarrow 0} \frac{-\langle \varphi,\ (f-t(J^*(\varphi_A)))^+ - f^+\rangle}{\|f-t(J^*(\varphi_A))-f\|_p + \|(f-t(J^*(\varphi_A)))^+ - f^+\|_p}$$

$$= \limsup_{t\downarrow 0} \frac{-\langle \varphi,\ f^+ - t(J^*(\varphi_A)) - f^+\rangle}{\|f-t(J^*(\varphi_A))-f\|_p + \|f^+ - t(J^*(\varphi_A)) - f^+\|_p}$$

$$= \limsup_{t\downarrow 0} \frac{t\langle \varphi, J^*(\varphi_A)\rangle}{2t\|J^*(\varphi_A)\|_p}$$

$$= \frac{\langle \varphi, J^*(\varphi_A)\rangle}{2\|J^*(\varphi_A)\|_p}$$

$$= \frac{\langle \varphi_A, J^*(\varphi_A)\rangle}{2\|J^*(\varphi_A)\|_p}$$

$$= \frac{\|J^*(\varphi_A)\|_p}{2}$$

$$> 0.$$

This contradiction proves ($B_1$). That is,

$$\theta^* \in \widehat{D}^* P_{K_p}(f)(\varphi) \implies \mu(\{s\in S: \varphi(s) < 0\} \cap \{s\in S: f(s) > 0\}) = 0.$$

Proof of ($B_2$). Assume, by contradiction that $\mu(\{s\in S: \varphi(s) > 0\} \cap \{s\in S: f(s) > 0\}) > 0$. Since $\varphi(s) > 0$ if and only if $J^*(\varphi)(s) > 0$. This implies

$$\mu(\{s\in S: J^*(\varphi)(s) > 0\} \cap \{s\in S: f(s) > 0\}) > 0.$$

Then, there are positive numbers $\alpha, \beta$ such that

$$\mu(\{s\in S: 0 < J^*(\varphi)(s) < \alpha\} \cap \{s\in S: f(s) > \beta\}) > 0.$$

Let $B := \{s\in S: 0 < J^*(\varphi)(s) < \alpha\} \cap \{s\in S: f(s) > \beta\}$. Define

$$\varphi_B(s) = \begin{cases} \varphi(s), & \text{if } s \in B, \\ 0, & \text{if } s \notin B. \end{cases}$$

In this case, we take a directional line segment in the limit (5.13) as, $h = f - tJ^*(\varphi_B)$, for $t \downarrow 0$ with $0 < t < \frac{\beta}{\alpha}$. Then, one has

$$(f - t(J^*(\varphi_B)))^+ = f^+ - t(J^*(\varphi_B)).$$

This implies

$$0 \geq \limsup_{h\to f} \frac{-\langle \varphi,\ h^+ - f^+\rangle}{\|h-f\|_p + \|h^+ - f^+\|_p}$$

$$\geq \limsup_{t\downarrow 0, t < \frac{\beta}{\alpha}} \frac{-\langle \varphi,\ (f-t(J^*(\varphi_B)))^+ - f^+\rangle}{\|f-t(J^*(\varphi_B))-f\|_p + \|(f-t(J^*(\varphi_B)))^+ - f^+\|_p}$$

$$= \limsup_{t\downarrow 0, t<\frac{\beta}{\alpha}} \frac{-\langle \varphi, f^+ - t(J^*(\varphi_B)) - f^+\rangle}{\|f - t(J^*(\varphi_B)) - f\|_p + \|f^+ - t(J^*(\varphi_B)) - f^+\|_p}$$

$$= \limsup_{t\downarrow 0, t<\frac{\beta}{\alpha}} \frac{t\langle \varphi, J^*(\varphi_B)\rangle}{2t\|J^*(\varphi_B)\|_p}$$

$$= \frac{\langle \varphi, J^*(\varphi_B)\rangle}{2\|J^*(\varphi_B)\|_p}$$

$$= \frac{\langle \varphi_A, J^*(\varphi_B)\rangle}{2\|J^*(\varphi_B)\|_p}$$

$$= \frac{\|J^*(\varphi_B)\|_p}{2}$$

$$> 0.$$

This contradiction proves (B$_2$). That is,

$$\theta^* \in \widehat{D}^* P_{K_p}(f)(\varphi) \implies \mu(\{s \in S: \varphi(s) > 0\} \cap \{s \in S: f(s) > 0\}) = 0.$$

Proof of (B$_3$). Assume that

$$\mu(\{s \in S: \varphi(s) < 0\} \cap \{s \in S: f(s) \leq 0\}) > 0.$$

This is equivalent to

$$\mu(\{s \in S: J^*(\varphi)(s) < 0\} \cap \{s \in S: f(s) \leq 0\}) > 0.$$

Then, there are positive numbers $\alpha, \beta$ such that

$$\mu(\{s \in S: J^*(\varphi)(s) \leq -\alpha\} \cap \{s \in S: -\beta \leq f(s) \leq 0\}) > 0.$$

Let $E := \{s \in S: J^*(\varphi)(s) \leq -\alpha\} \cap \{s \in S: -\beta \leq f(s) \leq 0\}$. We take a sequence $\{E_n\}$ of subsets of $E$ with positive measures and satisfying $\mu(E_n) \downarrow 0$. For any given positive integer $n$, we define

$$f_n(s) = \begin{cases} f(s) + 2\beta, & \text{if } s \in E_n, \\ 0, & \text{if } s \notin E_n. \end{cases}$$

Since $f^+(s) = 0$, for all $s \in E$, this implies

$$(f + f_n)^+(s) - f^+(s) = \begin{cases} f(s) + 2\beta, & \text{if } s \in E_n, \\ 0, & \text{if } s \notin E_n \end{cases}$$

$$= f_n(s), \quad \text{for } n = 1, 2, \ldots.$$

In this case, we take a sequence directional limit in (5.13) as $h = f + f_n$, for $n \to \infty$. One has

$$0 \geq \limsup_{h \to f} \frac{-\langle \varphi, h^+ - f^+ \rangle}{\|h-f\|_p + \|h^+ - f^+\|_p}$$

$$\geq \limsup_{n \to \infty} \frac{-\langle \varphi, (f+f_n)^+ - f^+ \rangle}{\|f+f_n-f\|_p + \|(f+f_n)^+ - f^+\|_p}$$

$$= \limsup_{n \to \infty} \frac{-\langle \varphi, f_n \rangle}{\|f_n\|_p + \|f_n\|_p}$$

$$= \limsup_{n \to \infty} \frac{-\int_S \varphi(s) f_n(s) \mu(ds)}{2\|f_n\|_p}$$

$$= \limsup_{n \to \infty} \frac{-\int_{E_n} \varphi(s) f_n(s) \mu(ds)}{2\|f_n\|_p}$$

$$= \limsup_{n \to \infty} \frac{-\int_{E_n} \varphi(s)(f(s) + 2\beta) \mu(ds)}{2\|f_n\|_p}$$

$$\geq \limsup_{n \to \infty} \frac{\int_{E_n} \alpha\beta \mu(ds)}{2\|f_n\|_p}$$

$$= \limsup_{n \to \infty} \frac{\alpha\beta\mu(E_n)}{2\|f_n\|_p}$$

$$\geq \limsup_{n \to \infty} \frac{\alpha\beta\mu(E_n)}{4\beta\mu(E_n)}$$

$$= \frac{\alpha}{4}$$

$$> 0.$$

This contradiction proves (B$_3$). That is,

$$\theta^* \in \widehat{D}^* P_{K_p}(f)(\varphi) \implies \mu(\{s \in S: \varphi(s) < 0\} \cap \{s \in S: f(s) \leq 0\}) = 0.$$

Hence, the part "$\implies$" of (b) is proved. Next, we secondly prove the part "$\impliedby$" of (b). Suppose that the following equations are satisfied.

$$\mu(\{s \in S: \varphi(s) \neq 0 \text{ and } f(s) > 0\}) = 0 \quad \text{and} \quad \mu(\{s \in S: \varphi(s) < 0 \text{ and } f(s) \leq 0\}) = 0.$$

The first equation of the above two equations implies

$$f(s) > 0 \implies \varphi(s) = 0, \text{ for } \mu\text{-almost all } s \in S.$$

Since, for any $h \in l_p$, we have $h^+(s) \geq 0$, for all $s \in S$, it follows

$$-\langle \varphi, h^+ - f^+ \rangle$$

$$= -\int_S \varphi(s)(h^+(s) - f^+(s)) \mu(ds)$$

$$= -\int_{f(s)>0} \varphi(s)(h^+(s) - f^+(s))\, \mu(ds) - \int_{f(s)\leq 0} \varphi(s)(h^+(s) - f^+(s))\mu(ds)$$

$$= 0 - \int_{f(s)\leq 0} \varphi(s)\big(h^+(s) - f^+(s)\big)\mu(ds)$$

$$= -\int_{f(s)\leq 0 \text{ and } \varphi(s)<0} \varphi(s)h^+(s)\mu(ds) - \int_{f(s)\leq 0 \text{ and } \varphi(s)\geq 0} \varphi(s)h^+(s)\mu(ds)$$

$$= 0 - \int_{f(s)\leq 0 \text{ and } \varphi(s)\geq 0} \varphi(s)h^+(s)\mu(ds)$$

$$\leq 0.$$

This implies

$$\limsup_{h\to f} \frac{-\langle \varphi,\ h^+ - f^+\rangle}{\|h-f\|_p + \|h^+ - f^+\|_p}$$

$$\leq \limsup_{t\downarrow 0} \frac{0}{\|h-f\|_p + \|h^+ - f^+\|_p}$$

$$= 0.$$

By (5.13), it follows that

$$\mu(\{s \in S: \varphi(s) \neq 0 \text{ and } f(s) > 0\}) = 0$$
$$\text{and } \mu(\{s \in S: \varphi(s) < 0 \text{ and } f(s) \leq 0\}) = 0$$
$$\Longrightarrow \quad \theta^* \in \widehat{D}^* P_{K_p}(f)(\varphi).$$

The proof of (b) is completed.

Proof of (c). For any $f \in K_p \setminus \{\theta\}$, we have

$$\mu(\{s \in S: J(f)(s) > 0 \text{ and } f(s) > 0\}) = \mu(\{s \in S: f(s) > 0\}) > 0.$$

Then, by (b), we have

$$\theta^* \notin \widehat{D}^* P_{K_p}(f)(J(f)), \text{ for any } f \in K_p \setminus \{\theta\}.$$

Let $f \in K_p$. For an arbitrarily given $\psi \in K_q$, and for $\varphi \in L_q(S)$, by (2.20), we have

$$\varphi \in \widehat{D}^* P_{K_p}(f)(\psi) \iff \limsup_{h\to f} \frac{\langle(\varphi,-\psi),\ (h, P_{K_p}(h)) - (f, P_{K_p}(f))\rangle}{\|h-f\|_p + \|P_{K_p}(h) - P_{K_p}(f)\|_p} \leq 0$$

$$\iff \limsup_{h\to f} \frac{\langle \varphi,\ h-f\rangle - \langle \psi,\ h^+ - f^+\rangle}{\|h-f\|_p + \|h^+ - f^+\|_p} \leq 0$$

$$\iff \limsup_{h\to f} \frac{\langle \varphi,(h^+-f^+)+(h^--f^-)\rangle - \langle \psi,\ h^+-f^+\rangle}{\|h-f\|_p + \|h^+-f^+\|_p} \leq 0$$

$$\Leftrightarrow \quad \limsup_{h\to f} \frac{\langle \varphi, h^- - f^-\rangle + \langle \varphi-\psi, h^+ - f^+\rangle}{\|h-f\|_p + \|h^+ - f^+\|_p} \leq 0. \qquad (5.14)$$

In particular,

$$\varphi \in \widehat{D}^* P_{K_p}(f)(J(f)) \quad \Leftrightarrow \quad \limsup_{h\to f} \frac{\langle \varphi, h^- - f^-\rangle + \langle \varphi-J(f), h^+ - f^+\rangle}{\|h-f\|_p + \|h^+ - f^+\|_p} \leq 0. \qquad (5.15)$$

Since $f \in K_p$, then, $f^- = \theta$ and $J(f) \in K_q$. Notice that $h^- \in -K_p$, for any $h \in K_p$. Substituting $\varphi$ by $J(f)$ in the limit in (5.15), we have

$$\limsup_{h\to f} \frac{\langle (J(f), -J(f)), (h, P_{K_p}(h)) - (f, P_{K_p}(f))\rangle}{\|h-f\|_p + \|P_{K_p}(h) - P_{K_p}(f)\|_p}$$

$$= \limsup_{h\to f} \frac{\langle J(f), h^- - f^-\rangle + \langle J(f) - J(f), h^+ - f^+\rangle}{\|h-f\|_p + \|h^+ - f^+\|_p}$$

$$= \limsup_{h\to f} \frac{\langle J(f), h^- - f^-\rangle}{\|h-f\|_p + \|h^+ - f^+\|_p}$$

$$= \limsup_{h\to f} \frac{\langle J(f), h^-\rangle}{\|h-f\|_p + \|h^+ - f^+\|_p}$$

$$= \limsup_{h\to f} \frac{\int_S J(f)(s) h^-(s) \mu(ds)}{\|h-f\|_p + \|h^+ - f^+\|_p}$$

$$\leq \limsup_{h\to f} \frac{0}{\|h-f\|_p + \|h^+ - f^+\|_p}$$

$$= 0.$$

By (5.14), we obtain

$$J(f) \in \widehat{D}^* P_{K_p}(f)(J(f)), \text{ for any } f \in K_p. \qquad (5.16)$$

Proof of (d). Let $\psi \in K_q$. For any $\varphi \in K_q$ satisfying $\varphi \leqslant_q \psi$, replacing $f$ by $\theta$ in the limit (5.14), we have

$$\varphi \in \widehat{D}^* P_{K_p}(\theta)(\psi) \quad \Leftrightarrow \quad \limsup_{h\to \theta} \frac{\langle \varphi, h^- - \theta^-\rangle + \langle \varphi-\psi, h^+ - \theta^+\rangle}{\|h-\theta\|_p + \|h^+ - \theta^+\|_p} \leq 0.$$

$$\Leftrightarrow \quad \limsup_{h\to \theta} \frac{\langle \varphi, h^-\rangle + \langle \varphi-\psi, h^+\rangle}{\|h\|_p + \|h^+\|_p} \leq 0. \qquad (5.17)$$

By $\varphi \in K_q$ satisfying $\varphi \leqslant_q \psi$, it follows that

$$\limsup_{h\to \theta} \frac{\langle \varphi, h^-\rangle + \langle \varphi-\psi, h^+\rangle}{\|h\|_p + \|h^+\|_p}$$

$$\leq \limsup_{h \to \theta} \frac{0+0}{\|h\|_p + \|h^+\|_p}$$

$$= 0.$$

For an arbitrarily given $\psi \in K_q$, by (5.17), this implies

$$\varphi \in \widehat{D}^* P_{K_p}(\theta)(\psi), \text{ for any } \varphi \in K_q \text{ with } \varphi \leqslant_q \psi.$$

This induces that, for any given $\psi \in K_q$, we have

$$[\theta^*, \psi]_{\leqslant_q} \subseteq \widehat{D}^* P_{K_p}(\theta)(\psi). \tag{5.18}$$

Next, we prove the opposite inclusion of (5.18). For this arbitrarily given $\psi \in K_q$, let $\varphi \in L_q$. Suppose $\varphi \notin [\theta^*, \psi]_{\leqslant_q}$. Then, rest of the proof of part (d) is divided to the following two cases.

Case 1. $\varphi \not\geqslant_q \theta^*$. In this case, $\mu\{s \in S : \varphi(s) < 0\} > 0$. Let $A := \{s \in S : \varphi(s) < 0\}$. Define

$$\varphi_A(s) = \begin{cases} \varphi(s), & \text{for } s \in A, \\ 0, & \text{for } s \notin A. \end{cases}$$

In this case, we take a directional line segment in the limit (5.17) as $h = tJ^*(\varphi_A)$, for $t \downarrow 0$. Then, one has

$$h^+ = \theta \text{ and } h^- = tJ^*(\varphi_A), \text{ for all } t > 0.$$

This implies

$$\limsup_{h \to \theta} \frac{\langle \varphi, h^- \rangle + \langle \varphi - \psi, h^+ \rangle}{\|h\|_p + \|h^+\|_p}$$

$$\geq \limsup_{t \to 0} \frac{\langle \varphi, tJ^*(\varphi_A) \rangle + \langle \varphi - \psi, \theta \rangle}{\|tJ^*(\varphi_A)\|_p + \|\theta\|_p}$$

$$= \limsup_{t \to 0} \frac{t \langle \varphi, J^*(\varphi_A) \rangle}{t\|J^*(\varphi_A)\|_p}$$

$$= \frac{\langle \varphi, J^*(\varphi_A) \rangle}{\|J^*(\varphi_A)\|_p}$$

$$= \frac{\langle \varphi_A, J^*(\varphi_A) \rangle}{\|J^*(\varphi_A)\|_p}$$

$$= \|J^*(\varphi_A)\|_p$$

$$> 0.$$

This contradiction shows that

$$\varphi \not\geqslant_q \theta^* \implies \varphi \notin \widehat{D}^* P_{K_p}(\theta)(\psi). \tag{5.19}$$

However, from case 1, we have that to get $\varphi \in \widehat{D}^* P_{K_p}(\theta)(\psi)$, it is necessary to have $\varphi \succcurlyeq_q \theta$. That is, from case 1, we have

$$\varphi \in \widehat{D}^* P_{K_p}(\theta)(\psi) \Longrightarrow \varphi \in K_q. \tag{5.20}$$

Case 2. $\varphi \not\preccurlyeq_q \psi$ (This case includes that $\varphi >_q \psi$ or $\varphi$ and $\psi$ are not $\preccurlyeq_q$-compareable). Then, by (5.20), in case 2, we have $\mu\{s \in S : \varphi(s) - \psi(s) > 0\} > 0$. Let $B := \{s \in S : \varphi(s) - \psi(s) > 0\}$. Then $\mu(B) > 0$. Since $\psi \in K_q$, it implies that

$$J^*(\varphi)(s) > 0, \text{ for all } s \in B. \tag{5.21}$$

Define

$$\varphi_B(s) = \begin{cases} \varphi(s), & \text{for } s \in B, \\ 0, & \text{for } s \notin B. \end{cases}$$

And

$$(\varphi - \psi)_B(s) = \begin{cases} (\varphi - \psi)(s), & \text{for } s \in B, \\ 0, & \text{for } s \notin B. \end{cases} \tag{5.22}$$

In this case, we take a directional line segment in the limit (5.17) as $h = tJ^*(\varphi_B)$, for $t \downarrow 0$. Then, one has

$$h^+ = tJ^*(\varphi_B) \quad \text{and} \quad h^- = \theta, \text{ for all } t > 0.$$

By (5.21) and (5.22), this implies

$$\limsup_{h \to \theta} \frac{\langle \varphi, h^- \rangle + \langle \varphi - \psi, h^+ \rangle}{\|h - \theta\|_p + \|h^+ - \theta^+\|_p}$$

$$= \limsup_{h \to f} \frac{\langle \varphi, \theta \rangle + \langle \varphi - \psi, tJ^*(\varphi_B) \rangle}{\|\theta\|_p + t\|J^*(\varphi_B)\|_p}$$

$$= \limsup_{h \to f} \frac{t\langle \varphi - \psi, J^*(\varphi_B) \rangle}{t\|J^*(\varphi_B)\|_p}$$

$$= \frac{\langle \varphi - \psi, J^*(\varphi_B) \rangle}{\|J^*(\varphi_B)\|_p}$$

$$= \frac{\langle (\varphi - \psi)_B, J^*(\varphi_B) \rangle}{\|J^*(\varphi_B)\|_p}$$

$$> 0.$$

This contradiction shows that

$$\varphi \not\preccurlyeq_q \psi \implies \varphi \notin \widehat{D}^* P_{K_p}(\theta)(\psi). \tag{5.23}$$

By (5.18), (5.19) and (5.23), part (d) is proved. □

## Acknowledgments

The author is very grateful to Professor Boris S. Mordukhovich for his kind communications, valuable suggestions and enthusiasm encouragements in the development stage of this paper.
## References

[1] Alber, Ya., Metric and generalized projection operators in Banach spaces: properties and applications, in "Theory and Applications of Nonlinear Operators of Accretive and Monotone Type" (A. Kartsatos, Ed.), Marcel Dekker, inc. (1996), 15–50.

[2] Alber, Ya. I., James's orthogonality and orthogonal decompositions of Banach spaces, J. Math. Anal. Appl. 312 (2005) 330–342.

[3] Bauschke H. H. and Combettes, P. L., Convex analysis and monotone operator theory in Hilbert spaces, ISSN 1613-5237 Springer, Heidelberg New York Dordrecht London (2011).

[4] Berdyshev, V.I., Differentiability of the metric projection in normed spaces. Collection: Approximation of functions by polynomials and splines, 150, 58−71. Akad. Nauk. SSSR, Ural. Nauchn. Tsentr., Sverd. (1985).

[5] Braess, D., Nonlinear Approximation Theory, Springer, Berlin, 1986.

[6] Fitzpatrick, S. and R.R. Phelps, Differentiability of the metric projection in Hilbert space, Trans. Amer. Math. Soc, 270 (1982), 483−501.

[7] Giles, J. R., On a characterization of differentiability of the norm of a normed linear space, (Received 5 March 1969 https://doi.org/10.1017/S1446788700008387 Published online by Cambridge University Press.

[8] Goebel, K. and Reich, S., Uniform Convexity, Hyperbolic Geometry, and Nonexpansive Mappings, Marcel Dekker, New York and Basel (1984).

[9] Haraux, A., How to differentiate the projection on a convex set in Hilbert space. Some applications to variational inequalities. J. Math. Soc. Japan, 29 (1977), 615−631

[10] Holmes, R. B., Smoothness of certain metric projections on Hilbert space, Trans. Amer. Math. Soc, 184 (1973), 87−100.

[11] Khan, A. A., Li, J. L. and Reich, S., Generalized projections on general Banach spaces, J. Nonlinear Convex Anal. 24 (2023), 1079—1112.

[12] Li, J. L., Directional Differentiability of the Metric Projection Operator in Uniformly Convex and Uniformly Smooth Banach Spaces, Journal of Optimization Theory and Applications.

[13] Li, J. L., Strict Fréchet differentiability of the metric projection operator in Hilbert spaces, submitted.

[14] Li, J. L., Mordukhovich derivatives of the metric projection operator in Hilbert spaces, submitted.

[15] Li, J. L., Fréchet differentiability of the metric projection operator in Banach spaces, submitted.

[16] Lindenstrauss, Joram; Tzafriri, Lior, Classical Banach spaces. II. Function spaces, Ergebnisse der Mathematik und ihrer Grenzgebiete [Results in Mathematics and Related Areas], vol. 97, Berlin-New York: Springer-Verlag, pp. x+243, ISBN 3-540-08888-1. (1979).

[17] Malanowski, K., Differentiability of projections onto cones and sensitivity analysis for optimal control, Proceedings of the 41st IEEE Conference on Decision and Control, Las Vegas, Nevada USA, December (2002).

[18] Mordukhovich, Boris S. Variational Analysis and Generalized Differentiation I, Basic Theory. DOI 10.1007/978-3-540-31247-5 Springer Heidelberg New York Dordrecht London (2006).

[19] Noll, D., Directional differentiability of the metric projection in Hilbert space, Pacific Journal of Mathematics, Vol. 170 No. 2 (1995).

[20] Reich, S., A remark on a problem of Asplund, Atti Accad. Lincei 67 (1979), 204—205,

[21] Shapiro, A., On differentiability of the metric projection in W1. Boundary case. Proc. Amer. Math. Soc, 99 (1987), 123−128.

[22] Shapiro, A., On concepts of directional differentiability, J. Math. Anal. Appl., 86 (1990), 77−487.

[23] Takahashi, W., Nonlinear Functional Analysis, Yokohama Publishers, (2000).